\documentclass[11pt]{article}
\usepackage{a4wide}                 
\usepackage{lmodern}                

\usepackage{amsmath}
\usepackage{amsthm}
\usepackage{amsfonts}
\usepackage{amssymb}
\usepackage{mathrsfs}               
\usepackage{bbm}               
\usepackage{array}
\usepackage{graphicx}
\usepackage{natbib}
\usepackage{needspace}
\usepackage{multirow}               
\usepackage{hyperref}
\usepackage{longtable}

\input{math.sty}

\allowdisplaybreaks[3]

\title{Quantile estimation for L\'evy measures}
\author{Mathias Trabs\footnote{E-mail address: trabs@math.hu-berlin.de}~\footnote{ The author thanks Markus Rei{\ss} for helpful comments. This research was partially supported by the Deutsche Forschungsgemeinschaft through the SFB 649 ``Economic Risk''.}}

\date{Humboldt-Universit\"at zu Berlin}

\begin{document}

\maketitle

\begin{abstract}
  Generalizing the concept of quantiles to the jump measure of a L\'evy process, the generalized quantiles $q_{\tau}^{\pm}>0$, for $\tau>0$, are given by the smallest values such that a jump larger than $q_{\tau}^{+}$ or a negative jump smaller than $-q_{\tau}^{-}$, respectively, is expected only once in $1/\tau$ time units. Nonparametric estimators of the generalized quantiles are constructed using either discrete observations of the process or using option prices in an exponential L\'evy model of asset prices. In both models minimax convergence rates are shown. Applying Lepski's approach, we derive adaptive quantile estimators. The performance of the estimation method is illustrated in simulations and with real data. 
\end{abstract}



\section{Introduction}

Whenever the modeling of random processes in biology, finance or physics requires to incorporate jumps, L\'evy processes are one of the building blocks under consideration. Consequently, their statistical analysis attracted much attention in the last decades. The estimation of the jump distribution, characterized by the L\'evy measure, is of particular interest. So far, only the jump density or linear functionals of it like the corresponding distribution function has been the aim estimation procedures. In the present work, we study the estimation of the (generalized) quantiles of the L\'evy measure that is the inverse of the distribution function. For a given intensity $\tau>0$ the quantile $q_\tau^+$ is the minimal jump height such that jumps larger than $q_\tau^+$ are expected only $\tau$ times in one time unit. Equivalently a jump larger than $q_\tau^+$ is expected only once in $1/\tau$ time units. $q_\tau^-$ is analogously defined for negative jumps, see Section~\ref{sec:genQuan} for the precise definitions and a discussion of possible applications. 

We will consider two different observation schemes: If we directly observe the L\'evy process at equidistant discrete time points $\Delta,2\Delta,\dots ,n\Delta$, the estimator relies on the increments $L_{\Delta k}-L_{\Delta(k-1)}$, for $k=1,\dots,n$, which are independent and identically distributed according to the law of $L_\Delta$. We will focus on low-frequency observations where $\Delta>0$ remains fixed while $n\to\infty$. 

The second observation scheme is motivated by an application in finance. Modeling an asset as exponential of a L\'evy process $L$, we use prices of put and call options to estimate the characteristics of the L\'evy process under the risk-neutral measure. This allows to estimate how shocks, in the sense of large jumps, are priced into the asset. Since the observed option prices are noisy, we face a statistical problem of regression type. The error analysis of the nonparametric estimators is similar to the case of low frequent direct observations.

In a high-frequency regime, i.e. observing $L_{k\Delta}$ with $\Delta\downarrow 0$, we almost see the jumps in the path and thus the estimation of the L\'evy measure is relatively straight forward, see \cite{aitSahaliaJacod:2012} for a review. As noticed by \cite{neumannReiss2009} in the low-frequency regime, i.e. $\Delta>0$ is fixed, the nonparametric estimation problem is more difficult because the number of jumps that occurred within an increment is not identifiable. Using however the L\'evy--Khintchine formula which links the characteristic function of the marginal distribution and the characteristic triplet, we can estimate the jump measure by a spectral approach. This idea was initiated by \cite{belomestnyReiss2006} and then studied further, see \cite{reiss2013} for an overview. The observation scheme of the option prices was first studied by \cite{contTankov:2004} as well as \cite{belomestnyReiss2006}. 

With the notable exception of \cite{belomestny2010}, who estimates the fractional order of a L\'evy process, only linear functionals of the jump measure are considered. Instead, we solve the substantially more demanding problem of estimating the nonlinear generalized quantiles in this nonlinear inverse problem. The conditions we impose on the models are very weak. We allow the whole spectrum of L\'evy processes, reaching from diffusions to pure jump processes with unbounded variation and the combination of both. The quantile estimators are very robust in the sense that they do not depend on the drift and volatility parameters, which may be surprising. In both above described observation schemes we derive convergence rates for the quantile estimators (Theorems~\ref{thm:rateQuantiles} and \ref{thm:quantileFin}). In view of the literature our rates appear to be minimax optimal. As side results convergence rates for density estimation and distribution function estimation are obtained which are novel due to the 
studied generality of L\'evy processes. 

The question of adaptive estimation methods for L\'evy processes in the high-frequency and low-frequency regime has  been only recently addressed by \cite{comteGenonCatalot2011} and \cite{kappus2012}, respectively, who apply model selection procedures. In the exponential L\'evy model we provide an adaptive method of Lepski-type that achieves the optimal rates (Theorem~\ref{thm:quantileFinAdapt}) with an additional $(\log \log n)$-payment for adaptivity that appears to be unavoidable. Since the exponential L\'evy model is related to the low-frequency regime let us compare our results to one by \cite{kappus2012} who has estimated linear functionals of the L\'evy density. Profiting from the regression structure, we can handle a broad class of processes while results in \cite{kappus2012} are restricted to L\'evy processes of bounded variation. The model selection approach leads to a finite sample oracle inequality where the constants, however, are not explicitly determined which might be problematic for 
applications. Although our error analysis is completely asymptotic, our data-driven method achieves very good results in finite sample situations, as illustrated with prices of options on the German DAX-index.  

In the next Section the generalized quantiles are introduced and their applications are discussed. In addition the basic estimation idea is outlined. Discrete observations of the L\'evy process and the option price model are considered in Sections~\ref{sec:quantileEst} and \ref{sec:FinancialQuantile}, respectively. In the latter model we construct an adaptive version of the estimator in Section~\ref{sec:lepski}. Simulations and a real data study based on DAX-options are given in Section~\ref{sec:sim}. All proofs are postponed to Section~\ref{sec:LevyProofs}.

\section{Generalized quantiles and estimation principle}\label{sec:genQuan}

Due to the L\'evy--It\^o decomposition, any L\'evy process $L=\{L_t:t\ge0\}$ can be represented as the sum of a deterministic drift determined by a parameter $\gamma\in\R$, a Brownian motion with volatility $\sigma^2>0$ and an independent jump component. The distribution of the jump sizes is described by the L\'evy measure $\nu$ which may have a singularity at zero, but satisfies $\int (x^2\wedge 1)\nu(\d x)<\infty$. The process $L$ is uniquely determined by this so-called characteristic triplet $(\sigma^2,\gamma,\nu)$. 

By definition $\nu(A)$, for any Borel set $A\in\mathscr{B}(\R)$, is the expected number of jumps per unit time whose size belongs
to $A$. Taking into account the possible singularity of $\nu$ at
zero, the generalized distribution function is defined by
\begin{equation}\label{eq:distfct}
  N(t):=\begin{cases}
  \nu((-\infty,t]), & \text{for }t<0,\\
  \nu([t,\infty)), & \text{for }t>0.
  \end{cases}
\end{equation}
Our aim is to estimate its inverse function which we call \emph{generalized
quantile function}. For a given level $\tau\in(0,\nu(\R_\pm))$ we introduce
\[
  q_\tau^-:=\sup\big\{t\ge0:\nu\big((-\infty,-t]\big)\ge\tau\big\}\quad\text{and}\quad
  q_\tau^+:=\sup\big\{t\ge0:\nu\big([t,\infty)\big)\ge\tau\big\}.
\]
Hence, $q_{\tau}^{+}$ (resp. $q_{\tau}^{-}$) is the largest value such that jumps larger than $q_\tau^+$ (resp. smaller than $-q_\tau^-$) have a least intensity $\tau$.

Especially the distribution of large jumps, corresponding to small values of $\tau$, are conveniently described using generalized quantiles. Let us discuss a few possible applications. For L\'evy processes with compound Poisson jump component, the quantiles of the L\'evy measure, which then are a finite measures, have very similar properties as quantiles of probability distributions. The total jump intensity of negative and positive jumps is given by $q_0^-$ and $q_0^+$, respectively, and the quantiles can be used as measures of location, scale and skewness of the jump distribution.  

These properties extend to infinite measures, too. For instance, the Bowley skewness for probability measures directly generalizes to $|q_\tau^--q_\tau^+|/(q_\tau^-+q_\tau^-)$ for $\tau>0$ quantifying symmetry between positive and negative jumps of a L\'evy processes. Another remarkable property is illustrated in the following example:
\begin{example}
  A common construction of L\'evy measures is the so-called exponential tilting where a L\'evy measure $\nu$ is multiplied by an exponential factor, i.e., $\tilde\nu_\lambda(\d x):=e^{-\lambda |x|}\nu(\d x)$ for $\lambda>0$. Supposing $\nu$ has a Lebesgue density which is bounded outside of a neighborhood of the origin and which converges polynomially fast to zero as $|x|\to\infty$, it is easy to see that the quantiles $q^\pm_{\tau,\lambda}$ corresponding to $\tilde\nu_\lambda$ grow like $|\log \tau|/\lambda$ as $\tau\to0$. Ploting $q^\pm_{\tau,\lambda}$ against $q^\pm_{\tau,0}$, the parameter $\lambda$ thus approximately shows up as inverse of the slope for large values of $q_{\tau,0}$.
\end{example}

Therefore, estimators for the generalized quantiles are important descriptive statistics for the jump behavior of L\'evy processes. But there are further applications: To construct and estimate L\'evy copulas, the generalized quantiles are necessary, cf. \cite[Chap. 5]{contTankov:2004b} and \cite{bucherVetter2013}, and in the context of modeling prices processes, they can be used to estimate the risk within the model. One of the most popular risk measures is the value-at-risk at some level $\tau\in(0,1)$ which is given by the $(1-\tau)$-quantile of the distribution of the loss of the asset. The generalized quantiles of $\nu$ are a closely related concept which takes only the influence of shocks into account. Following this idea, the quantiles may also be useful for dynamic quantile hedging in the spirit of \cite{follmerLeukert1999}.

Now, how to estimate these quantiles? Before we rigorously introduce the estimators and study their asymptotic properties in the following two sections let us outline the general estimation principle. We follow a similar strategy as \citet{dattnerEtAl2013} who study quantile estimation in the classical deconvolution model. Compared to deconvolution, the L\'evy model is harder for two reasons. First, it is a nonlinear inverse problem such that we have to linearize the estimation error and the remainder needs extra care. Second, the underlying deconvolution problem is determined by the distribution of the process itself. Consequently, there is a strong interplay between the jump measure that we want to estimate and the underlying deconvolution operator.   

Assuming $\int x^{2}\nu(\d x)<\infty$, $\phi_t$ is given by the L\'evy--Khintchine representation in Kolmogorov's version:
\begin{equation}\label{eqLevyKhint}
\phi_{t}(u):=\E[e^{iuY_{t}}]=e^{t\psi(u)}\quad\text{with}\quad\psi(u):=-\frac{\sigma^{2}}{2}u^{2}+i\gamma u+\int\big(e^{iux}-1-iux\big)\nu(\d x).
\end{equation}
Differentiating twice the characteristic exponent $\psi$, we obtain the estimating equation
\begin{equation}\label{eq:psiPP}
\psi''(u)=-\sigma^{2}-\F[x^{2}\nu](u)=\frac{\phi_{t}''(u)\phi_{t}(u)-\phi'_t(u)^2}{t\phi_t^{2}(u)}.
\end{equation}
As starting point we need an estimator of the characteristic function $\phi_t$ of the marginal distribution $L_t$ of the L\'evy process for some $t>0$. Using discrete observations of the L\'evy process, $\phi_t$ can be estimated by the empirical measure of the increments. In the financial model the characteristic function can be estimated via the pricing formula that links $\phi_t$ and the option prices. To estimate $q_\tau^\pm$, we will then apply the following program:
\begin{enumerate}
 \item A density estimator for the jump measure $\nu$ can be constructed by replacing $\phi_t$ in \eqref{eq:psiPP} with its estimator, regularizing with a band limited kernel and applying the inverse Fourier transform to the left and right-hand side of \eqref{eq:psiPP}.
 \item A plug-in approach yields an estimator for the distribution function.
 \item The generalized $\tau$-quantiles can be estimated by minimizing the distance between the value of distribution function estimator and $\tau$. 
\end{enumerate}
The estimator for the jump density is similar to the estimators proposed in \citet{nicklReiss2012} as well as \cite{kappus2012}. However, both articles are restricted to L\'evy processes with bounded variation and thus use only the first derivative of $\psi$. Focusing on settings that allow for parametric rates, the distribution function estimation have been considered by \cite{nicklReiss2012} and \cite{nicklEtAl2013} in the low and high-frequency regime, respectively.

\section{Discrete observations of the process}\label{sec:quantileEst}

We observe $n\in\mathbb{N}$ increments of the L\'evy process $L$
at equidistant time points with observation distance $\Delta>0$:
\begin{align*}
Y_{k}: & =L_{\Delta k}-L_{\Delta(k-1)},\quad k=1,\dots,n.
\end{align*}
The law of $Y_{k}$ will be denoted by $P_{\Delta}$. Using the empirical characteristic function $\phi_{\Delta,n}(u)=\frac{1}{n}\sum_{k=1}^{n}e^{iuY_{k}}$, we obtain an empirical version of $\psi''$ from \eqref{eq:psiPP}:
\begin{equation*}
  \hat{\psi}_{n}''(u)=\frac{\phi_{\Delta,n}''(u)\phi_{t}(u)-\phi'_{\Delta,n}(u)^2}{t\phi_{\Delta,n}^{2}(u)}\1_{\{|\phi_{\Delta,n}(u)|\ge(\Delta n)^{-1/2}\}},
\end{equation*}
where we multiply with the indicator function to stabilize against large stochastic errors. We define the density estimator as
\[
\hat{\nu}_{h}(t):=-t^{-2}\F^{-1}\Big[\hat{\psi}_{n}''(u)\F K(hu)\Big](t),\qquad t\neq0,
\]
where $K$ is a band-limited kernel with bandwidth $h>0$ satisfying for some order $p\in\mathbb{N}$
\begin{gather}
\begin{split}
\int_{\R}K(x)\d x=1,\qquad\int x^{l}K(x)\d x=0\quad\text{ for }l=1,\dots,p,\\
\supp\F K\subset[-1,1],\quad x^{p+1}K(x)\in L^{1}(\R).\qquad
\end{split}
\label{eq:propKernel}
\end{gather}
Note that $\hat\nu_h$ depends neither on the unknown volatility $\sigma^{2}$ nor on the drift parameter $\gamma$. The distribution function can be estimated via the left and the right tail integrals
\begin{align}
\hat{N}_{h}(t) & =-\int g_{t}(x)\F^{-1}\Big[\hat{\psi}_{n}''(u)\F K(hu)\Big](x)\d x\quad\text{with\quad}g_{t}(x):=\begin{cases}
x^{-2}\1_{(-\infty,t]}, & t<0,\\
x^{-2}\1_{[t,\infty)}, & t>0.
\end{cases}\label{eq:DefDistEst}
\end{align}
Owing to $g_t\in L^1(\R)$, the estimator $\hat N_h$ is always well defined. 

Assuming absolute continuity of $\nu$ at $\pm q_{\tau}^{\pm}$
with respect to the Lebesgue measure, the generalized quantiles $q_{\tau}^{+}>0$ and $q_{\tau}^{-}>0$ for a given level $\tau\in(0,\nu(\R_\pm))$ are determined by
\[
N(-q_{\tau}^{-})=\tau=N(q_{\tau}^{+}).
\]
If $\nu$ has finite mass on the negative or the positive halfline,
the $\tau$-quantiles only exist if $\tau\le\nu(\R_-)$ or $\tau\le\nu(\R_+)$, respectively. On the other hand for processes with high jump activity, i.e., $x^2\nu(x)\to\infty$ as $|x|\to0$, the bias of the distribution function estimator $\hat N_h(t)$ explodes as $|t|\to0$, cf. Proposition~\ref{prop:bias}. Since $\nu$ is not known, it is thus reasonable to estimate $q_{\tau}^{\pm}\vee\eta_{n}$
for some threshold $\eta_{n}>0$ and any $\tau>0$ instead of the quantiles themselves in order to stabilize the estimation problem.
For finite jump activity processes with regular jump densities at zero, the threshold value $\eta_{n}$ may converge slowly
to zero as $n\to\infty$ and we conclude that $\tau>\nu(\R^\pm)$ holds almost surely if $\eta_n\to0$ and the estimators take the value $\eta_{n}$ for all $n$. 

Using $\hat N_h$, a minimum contrast estimation approach yields the quantile estimator
\[
\hat{q}_{\tau,h}^{\pm}:=\argmin_{t\in[\eta_{n},\infty)}|\hat{N}_{h}(\pm t)-\tau|
\]
for threshold values which are either fixed or which logarithmically decay to zero. Note that we only consider the distribution function outside a neighborhood of the origin and thus our estimator does not depend on the diffusion component which corresponds to a Dirac measure at zero with mass $\sigma^2$.

Since $\hat q_{\tau,h}$ is a point estimator, for any fixed $\tau$, nonparametric convergence rates depend on the local smoothness of $\nu$ and it is natural to assume H\"older regularity. Before specifying the exact nonparametric classes of L\'evy triplets that we will consider, let us introduce some notation. For an open set $U\subset\R$ the space of all functions continuous on $U$ is denoted by $C(U)$. The set of all functions which are H\"older-regular with index $s>0$ on $U$ is given by
\begin{align*}
  C^s(U)&:=\{f\in C(U):\|f\|_{C^s(U)}<\infty\}\quad\text{with}\quad\\
  \|f\|_{C^s(U)}&:=\sum_{k=0}^{\lfloor s\rfloor}\sup_{x\in U}|f^{(k)}(x)|+\sup_{x,y\in U:x\neq y}\frac{|f^{\lfloor s\rfloor}(x)-f^{(\lfloor s\rfloor)}(y)|}{|x-y|^{s-\lfloor s\rfloor}},
\end{align*}
where $\lfloor s\rfloor$ denotes the smallest integer strictly smaller than $s$. If not specified differently, $\|\cdot\|_{L^p}$ with $p\ge1$ denotes the $L^p$-norm on the whole real line with the usual notation $\|\cdot\|_\infty$ for the supremum norm. The bounded variation norm of a function $f$ on an interval $[a,b]$ for $a<b$ is defined as
\[
  \|f\|_{BV([a,b])}:=\sup\Big\{\sum_{i=1}^n|f(x_i)-f(x_{i-1})|:n\in\N,a\le x_1<\dots<x_n\le b\Big\}.
\]
Define for the open set $U$, regularity $s>0$, number of moments $m>0$ and radius $R>0$ the class of L\'evy triplets
\begin{align*}
  \mathcal{C}^{s}(m,U,R):= & \Big\{(\sigma^{2},\gamma,\nu)\Big|\sigma^{2}\in[0,R],\gamma\in\R,\|x^{m}\nu\|_{L^{1}}\le R,\notag\\
 & \qquad\qquad\quad\nu\text{ has a Lebesgue density on }U\text{ with }\|\nu\|_{C^{s}(U)}\le R\Big\}.
\end{align*}
For the error analysis we distinguish between infinitely divisible distributions with polynomially decaying characteristic functions and exponentially decaying characteristic functions. For $\alpha,\delta,r>0$ and $\beta\in(0,2]$ set
\begin{align}
\mathcal{D}^{s}(\alpha,m,U,R):= & \Big\{(0,\gamma,\nu)\in\mathcal{C}^{s}(m,U,R)\Big|\|(1+|\bull|)^{-\Delta\alpha}/\phi_{\Delta}\|_{\infty}\le R,\|x\nu\|_{\infty}\le R,\notag\\
 & \qquad\qquad\qquad\qquad\qquad\quad \|x\nu\|_{BV([-\delta,\delta])}\le R\Big\},\notag\allowdisplaybreaks\\
\mathcal{E}^{s}(\beta,m,U,r,R):= & \Big\{(\sigma^{2},\gamma,\nu)\in\mathcal{C}^{s}(m,U,R)\Big|\|\exp(-r\Delta|\bull|^{\beta})/\phi_{\Delta}\|_{\infty}\le R\Big\}.\label{eq:ClassDist}
\end{align}
We will see that the class $\mathcal{D}^{s}(\alpha,m,U,R)$ corresponds to mildly ill-posed estimation problems. As shown in \cite{trabs2014} a polynomial decay of $\phi_\Delta$ is, up to a mild regularity assumption, equivalent to $x\nu$ beeing of bounded variation near the origin and, of course, $\sigma=0$. Hence, the conditions imposed in $\mathcal D^s$ are quite natural. Moreover, they allow to apply the Fourier multiplier theorem in \cite{trabs2014}. The estimation problem in the class $\mathcal{E}^{s}(\beta,m,U,R)$ is severely ill-posed leading to logarithmic convergence rates. Because the rates are so slow, we only need very mild assumptions in $\mathcal{E}^{s}(\beta,m,U,r,R)$ without any additional condition on the jump measure. This class especially contains all infinite variation processes, noting that L\'evy processes with a diffusion component satisfy the decay condition on $\phi_\Delta$ only for $\beta=2$.
 
The convergence rates will be established in the $\mathcal O_P$-sense which corresponds to the loss function of confidence intervals. Since the results hold uniformly over the given classes of L\'evy processes, we define the \emph{uniform stochastic Landau symbol} $\mathcal O_{P,\Theta}$ over a parameter set $\Theta$: For random variables $(A_n)_{n\in\N}$ we write $A_n=\mathcal O_{P,\Theta}(1)$ if 
\begin{equation}\label{eq:UniformOp}
  \lim_{R\to\infty}\limsup_{n\to\infty}\sup_{\theta\in\Theta}P_\theta(A_n>R)=0.                                                                                                                                                                                                                                                                                                       
\end{equation}

Since our estimation procedure relies on a plug-in approach using the density estimator $\hat \nu_h$, we start with its asymptotic behavior beeing of independent interest. As in the analysis of the quantile estimator in the deconvolution model by \cite{dattnerEtAl2013}, the following proposition is the first building block for showing the rates of $\hat q_{\tau,h}$. We will need the result for the uniform loss.  
\begin{prop}\label{prop:DensLevyRate} 
  Let $\alpha,\beta,s,r,R>0, m>4$ and let the kernel satisfy (\ref{eq:propKernel})
  with order $p\ge s$ and let $U\subset\R$ be a bounded, open set
  which is bounded away from zero. Then we have for $n\to\infty$
  \begin{enumerate}
  \item  uniformly in $(\sigma^{2},\gamma,\nu)\in\mathcal{D}^{s}(\alpha,m,U,R)$ for $h=h_{n,\Delta}=(\frac{\log n\Delta}{n\Delta})^{1/(2s+2\Delta\alpha+1)}$ 
  \[
  \sup_{t\in U}|\hat{\nu}_{h}(t)-\nu(t)|=\mathcal{O}_{P,\mathcal{D}^{s}}\Big(\Big(\frac{\log n\Delta}{n\Delta}\Big)^{s/(2s+2\Delta\alpha+1)}\Big),
  \]

  \item uniformly in $(\sigma^{2},\gamma,\nu)\in\mathcal{E}^{s}(\beta,m,U,r,R)$ for $h=h_{n,\Delta}=(\frac{\delta\log n\Delta}{2r\Delta})^{-1/\beta}$ with $\delta\in(0,3/2)$
  \[
  \sup_{t\in U}|\hat{\nu}_{h}(t)-\nu(t)|
  =\mathcal{O}_{P,\mathcal{E}^{s}}\Big(\Big(\frac{\log n\Delta}{\Delta}\Big)^{-s/\beta}\Big).
  \]
  \end{enumerate}
\end{prop}
In the mildly ill-posed case the rates correspond to the deconvolution problem with an error distribution whose characteristic function decays with polynomial rate $\Delta\alpha$. It can easily be verified that the for the pointwise loss we have the same rates without the logarithm in (i). They coincide with the convergence rates for the pointwise loss by \citet[Thm. 3.5]{kappus2012}, who has considered only finite variation L\'evy processes. \cite{kappus2012diss} have shown that these rates are minimax optimal. In the classical density estimation the logarithm is known to be unavoidable for the uniform loss. In comparison to the deconvolution model, the logarithmic rates in (ii) appear to be sharp, too.

The distribution function estimator $\hat N_h$ was studied by \citet{nicklEtAl2013} in a high-frequency regime. For low frequency observations a modification of $\hat{N}_{h}$ was considered by \citet{nicklReiss2012}. In both articles a uniform central limit theorem has been established. To this end, assumptions have been imposed which ensure that the parametric rate can be attained. Therefore, it is of interest to derive convergence rates for $\hat{N}_{h}$ in more general situations. We remark that the following result could be strengthened to the uniform loss on $\R\setminus[-\eta,\eta]$ for any $\eta>0$, cf. \eqref{eq:ConsistencyBoundLevy} below.
\begin{prop}\label{prop:DistFunct}
  Let $U\subset\R$ be an open set, $t\in U$
  and $\alpha,\beta,s,r,R>0$ and $m>4$. Suppose the kernel satisfies
  (\ref{eq:propKernel}) with order $p\ge s+1$. Let $\Delta\in(0,1)$
  and $n\to\infty$. Then:
  \begin{enumerate}
  \item We obtain uniformly in $(\sigma^{2},\gamma,\nu)\in\mathcal{D}^{s}(\alpha,m,U,R)$ with the bandwidth $h=h_{n,\Delta}=(n\Delta)^{-1/(2s+(2\Delta\alpha\vee1)+1)}$ 
  \begin{align*}
  |\hat{N}_{h}(t)-N(t)| & =\mathcal{O}_{P,\mathcal{D}^{s}}(r_{n,\Delta}),\quad r_{n,\Delta}:=\begin{cases}
  (n\Delta)^{-(s+1)/(2s+2\Delta\alpha+1)}, & \text{for }\Delta\alpha>1/2,\\
  (n\Delta)^{-1/2}(\log n\Delta)^{1/2}, & \text{for }\Delta\alpha=1/2,\\
  (n\Delta)^{-1/2}, & \text{for }\Delta\alpha\in(0,1/2).
  \end{cases}
  \end{align*}
  \item The choice $h=h_{n,\Delta}=(\frac{\delta}{2r})^{-1/\beta}(\frac{\log(n\Delta)}{\Delta})^{-1/\beta}$,
  for any $\delta\in(0,3/2)$, yields uniformly in $(\sigma^{2},\gamma,\nu)\in\mathcal{E}^{s}(\beta,m,U,r,R)$
  \[
  |\hat{N}_{h}(t)-N(t)|
  =\mathcal{O}_{P,\mathcal{E}^{s}}\Big(\Big(\frac{\log(n\Delta)}{\Delta}\Big)^{-(s+1)/\beta}\Big).
  \]
  \end{enumerate}
\end{prop}
As expected, the convergence rates for $\hat N_n$ are faster than for density estimation because we gain one degree of smoothness. In particular, we achieve the parametric rate for a L\'evy process with very slowly decaying characteristic function or for $\Delta$ sufficiently small.

Following the standard M-estimation strategy, we use
a Taylor expansion to analyze the estimation error of the quantile estimators. Using that $\hat{N}_{h}'(t)=-\sign(t)\hat{\nu}_{h}(t)$ for $t\neq0$, we obtain 
\begin{align*}
  0\approx\hat{N}_{h}(\hat{q}_{\tau,h}^{+})-\tau=\hat{N}_{h}(q_{\tau}^{+})-N(q_{\tau}^{+})-(\hat{q}_{\tau,h}^{+}-q_{\tau}^{+})\hat{\nu}_{h}(\xi^+)
\end{align*}
for some intermediate point $\xi^+$ between $q_{\tau}^{+}$and $\hat{q}_{\tau,h}^{+}$
and similarly for $\hat{q}_{\tau,h}^{-}$. By continuity of $\hat N_h$ and the construction of $\hat q_{\tau,h}^\pm$ the probability of the event $\{\hat{N}_{h}(\hat{q}_{\tau,h}^{+})-\tau=0\}$ converges to one, cf. \eqref{eq:ZEst}. On this event the estimation error can therefore be represented as 
\begin{equation}
|\hat{q}_{\tau,h}^{\pm}-q_{\tau}  ^{\pm}|=\Big|\frac{\hat{N}_{h}(\pm q_{\tau}^{\pm})-N(\pm q_{\tau}^{\pm})+\tau-\hat{N}_{h}(\hat{q}_{\tau,h}^{+})}{\hat{\nu}_{h}(\xi^{\pm})}\Big|=\frac{|\hat{N}_{h}(\pm q_{\tau}^{\pm})-N(\pm q_{\tau}^{\pm})|}{|\hat{\nu}_{h}(\xi^{\pm})|},\label{eq:errorRep}
\end{equation}
for intermediate points $\xi^\pm$. If the denominator does not explode, $q_{\tau}^{\pm}$ can be estimated with the same rate as $N(\pm q_\tau^\pm)$. Since the convergence rates for the density estimator $\hat\nu_h$ and for the distribution function estimator $\hat{N}_{h}$ are already established, it remains to show consistency of the quantile estimator $\hat{q}_{\tau,h}^{\pm}$ itself. To this end, a minimal global regularity of $\nu$ is required. Note that assuming a bounded density would be far to restrictive for L\'evy measures. 

We need to specify our nonparametric classes further such that the
quantiles exist. Writing $\nu\in C^{s'}(\R\setminus[-\eta,\eta])$ for $s'\in(-1,0],\eta>0$ if the antiderivatives $y\mapsto\int_{-\infty}^{-y} \nu(\d x)$ and $y\mapsto\int^\infty_{y} \nu(\d x)$ are in $C^{1+s'}((\eta,\infty))$, we define for a given $\tau>0$ and with $\zeta>0$
\begin{align}
&\mathcal{\tilde{D}}_{\tau}^{s,s'}(\alpha,m,\zeta,\eta,R)\label{eq:Dtilde}\\
&\quad:= \Big\{(\sigma^{2},\gamma,\nu)\in\mathcal{D}^{s}\big(\alpha,m,(q_{\tau}^{+}-\zeta,q_{\tau}^{+}+\zeta)\cup(q_{\tau}^{-}-\zeta,q_{\tau}^{-}+\zeta),R\big)\Big|\notag\\
&\qquad\qquad\exists q_{\tau}^{+},q_{\tau}^{-}\in(\eta,\infty): N(-q_{\tau}^{-}) =\tau=N(q_{\tau}^{+}), \|\nu\|_{C^{s'}(\R\setminus[-\eta,\eta])}<R,  \nu(q_{\tau}^{\pm})>\tfrac{1}{R}\Big\},\notag\\
&\mathcal{\tilde{E}}_{\tau}^{s,s'}(\beta,m,\zeta,\eta,r,R)\label{eq:ClassQuantile}\\
&\quad:= \Big\{(\sigma^{2},\gamma,\nu)\in\mathcal{E}^{s}\big(\beta,m,(q_{\tau}^{+}-\zeta,q_{\tau}^{+}+\zeta)\cup(q_{\tau}^{-}-\zeta,q_{\tau}^{-}+\zeta),r,R\big) \Big|\notag\\
&\qquad\qquad\exists q_{\tau}^{+},q_{\tau}^{-}\in(\eta,\infty):N(-q_{\tau}^{-})=\tau=N(q_{\tau}^{+}), \|\nu\|_{C^{s'}(\R\setminus[-\eta,\eta])}<R, \nu(q_{\tau}^{\pm})>\tfrac{1}{R}\Big\}.\notag
\end{align}
The condition $\nu(q_{\tau}^{\pm})>\tfrac{1}{R}$ especially implies that the quantiles are unique. As expected from the representation \eqref{eq:errorRep}, we obtain the same rates for quantile estimation as for distribution function estimation. 
\begin{theorem}\label{thm:rateQuantiles}
  Let $\tau>0$ and $\alpha,\beta,s,\zeta,r,R>0,s'\in(-1,0]$
  and $m>4$. Suppose the kernel satisfies (\ref{eq:propKernel}) with
  order $p\ge s+1$ and $\eta_{n}>0$ with $\eta_{n}^{-1}\lesssim\log n$.
  Then we obtain for $\Delta\in(0,1)$ and $n\to\infty$:
  \begin{enumerate}
  \item The bandwidth $h=h_{n,\Delta}=(n\Delta)^{-1/(2s+(2\Delta\alpha\vee1)+1)}$ yields uniformly in $(\sigma^{2},\gamma,\nu)\in\mathcal{\tilde{D}}_{\tau}^{s,s'}(\alpha,m,\zeta,\eta_{n},R)$
  \begin{align*}
  |\hat{q}_{\tau,h}^{\pm}-q_{\tau}^{\pm}| & =\mathcal{O}_{P,\mathcal{\tilde{D}}_{\tau}^{s,s'}}(r_{n,\Delta}),\quad r_{n,\Delta}:=\begin{cases}
  (n\Delta)^{-(s+1)/(2s+2\Delta\alpha+1)} & \text{for }\Delta\alpha>1/2,\\
  (n\Delta)^{-1/2}(\log n\Delta)^{1/2} & \text{for }\Delta\alpha=1/2,\\
  (n\Delta)^{-1/2} & \text{for }\Delta\alpha\in(0,1/2).
  \end{cases}
  \end{align*}
  \item The choice $h=h_{n,\Delta}=(\frac{\delta}{2r})^{-1/\beta}(\frac{\log(n\Delta)}{\Delta})^{-1/\beta}$,
  for any $\delta\in(0,2/3)$, yields uniformly in $(\sigma^{2},\gamma,\nu)\in\mathcal{\tilde{E}}_{\tau}^{s,s'}(\beta,m,\zeta,\eta_{n},r,R)$
  \[
  |\hat{q}_{\tau,h}^{\pm}-q_{\tau}^{\pm}|=\mathcal{O}_{P,\mathcal{\tilde{E}}_{\tau}^{s,s'}}\Big(\Big(\frac{\log(n\Delta)}{\Delta}\Big)^{-(s+1)/\beta}\Big).
  \]
  \end{enumerate}
\end{theorem}

\begin{remark}
  In the parametric regime., i.e., $(0,\gamma,\nu)\in\mathcal{\tilde{D}}_{\tau}^{s,s'}(\alpha,m,\zeta,\eta_{n},R)$ with $\Delta\alpha\in(0,1/2)$, we could hope for more. In view of the central limit theorem by \cite{nicklReiss2012} a more precise analysis of the main stochastic error term, defined in \eqref{eq:DecompLin}, may lead to a central limit theorem of the distribution function estimator $\hat N_{h}$. Since the our analysis of \eqref{eq:errorRep} reveals that
  \[
    \hat{q}_{\tau,h}^{+}-q_{\tau}^+=\frac{\hat{N}_{h}(q_{\tau}^{+})-N(q_{\tau}^{+})+o_P(n^{-1/2})}{\nu(q_\tau^+)+o_P(1)},
  \]
  a central limit theorem for $\hat q_\tau^\pm$ immediately follows. Note that optimality of the asymptotic variance by \cite{nicklReiss2012} in the sense of semi-parametric efficiency is proved in \cite{trabs2013}. Using the $\Delta$-method, this lower bound can be extend to a lower bound for the asymptotic variance for the quantile estimation.
\end{remark}

\begin{remark}
For high-frequency data, that is for $\Delta\to0$ and $n\Delta\to\infty$, we obtain in $\mathcal{\tilde{D}}_{\tau}^{s,s'}(\alpha,\zeta,\eta_{n},R)$
always the parametric rate $(n\Delta)^{-1/2}$. In $\mathcal{\tilde{E}}_{\tau}^{s,s'}(\beta,\zeta,\eta_{n},r,R)$ the bandwidth
$h_{n,\Delta}=\Delta^{1/\beta\wedge 3/4}$ should be chosen instead such that our estimates in the proof of Theorem~\ref{thm:rateQuantiles} yield the almost optimal rate $(n\Delta)^{-1/2}|\log\Delta|$ provided that
$s$ is large enough such that the bias condition $h^{s+1}=\mathcal{O}((n\Delta)^{-1/2})$ is satisfied. With stronger assumptions the latter restriction can be circumvented and the rate can be improved. In fact, in \citet{nicklEtAl2013} we show that the parametric rate $(n\Delta)^{-1/2}$ can be obtained with this estimator under suitable conditions on $\nu$. 
\end{remark}

\section{Observation of option prices}\label{sec:FinancialQuantile}

Let us consider the exponential L\'evy model for asset prices
\begin{equation}
S_{t}=S_{0}e^{rt+L_{t}},\quad t\ge0,\label{eq:expLevyMod}
\end{equation}
with initial value $S_{0}>0$, riskless interest rate $r\ge0$ and with the driving L\'evy process $L$ whose characteristic triplet is $(\sigma^{2},\gamma,\nu)$. Since L\'evy processes are a quite large and flexible class of stochastic processes, important stylized facts of financial data can be reproduced, see \cite{contTankov:2004b} for properties and examples of the exponential L\'evy model. At the same time the well understood probabilistic structure of L\'evy processes allows to construct procedures to fit the model to real data.  The nonparametric calibration of this model was studied by \cite{contTankov:2004,belomestnyReiss2006} as well as \cite{trabs:2011}. \cite{soehl2014} have derived confidence sets. An empirical study of the calibration methods can be found in \cite{soehlTrabs2014}. 

By arbitrage arguments the price process of an asset should be a martingale implying $\E[S_t]=S_0$ for all $t>0$. Assuming again $\int x^2\nu(\d x)<\infty$, this is equivalent to the martingale condition
\begin{equation}\label{eqMartCond}
    \frac{\sigma^2}{2}+\gamma+\int_{-\infty}^\infty(e^x-1-x)\nu(\d x)=0,
\end{equation}
which we assume throughout this section. 
Since we want to estimate L\'evy measure under the risk-neutral measure, the procedure is based on option prices. More precisely, we observe prices of vanilla options. Let us fix a maturity $T>0$, measured in years, and define the negative log--moneyness $x:=\log(K/S_0)-rT$ as a logarithmic transform of the strike prices $K>0$. In terms of $x$, call and put prices in model \eqref{eq:expLevyMod} are given by $\mathcal{C}(x,T)=S_0\mathbb E[(e^{L_T}-e^x)_+]$ and $\mathcal{P}(x,T)=S_0\mathbb E[(e^x-e^{L_T})_+]$, respectively. The prices can be summarized in the \emph{option function}
\begin{equation}\label{eq:optionFct}
  O(x):=\begin{cases}
                     \displaystyle S_0^{-1}\mathcal{C}(x,T), & \quad x\ge0,\\
                     \displaystyle S_0^{-1}\mathcal{P}(x,T), & \quad x<0.
                  \end{cases}
\end{equation}
We observe $O$ at a finite number of (transformed) strike prices $x_1,\dots,x_n$, for $n\in\N$, corrupted by noise 
\begin{equation}\label{eq:Observations}
  O_{j}=O(x_{j})+\delta_{j}\eps_j\quad j=1,\dots,n,
\end{equation}
where $(\eps_{j})$ are i.i.d. centered random variables with $\Var(\eps_j)=1$ and with local
noise levels $(\delta_j)$. The observation errors are due to the bid--ask spread and other market frictions. By interpolating the observations $(x_{j},O_{j})_{j=1,\dots,n}$ using B-splines, we construct an empirical version $\tilde{O}$ of the option function as in \cite{belomestnyReiss2006}. From the theoretical perspective linear splines are sufficient, but in applications B-splines of higher degrees may lead to better results, cf. \cite{soehlTrabs2014} for details. Since the function $O$ is related to the characteristic function $\phi_T$ of $L_{T}$ via the pricing formula, cf. \cite{carrMadan:1999},
\begin{equation}\label{pricingFormula}
  \mathcal{F}O(u):=\int_{-\infty}^{\infty}e^{iux}O(x)\d x=\frac{1-\phi_T(u-i)}{u(u-i)},
\end{equation} 
we obtain an estimator of $\phi_T$ given by 
\begin{equation}\label{eq:phiTilde}
\tilde\phi_{T,n}(u):=1-u(u+i)\F \tilde O(u+i).
\end{equation}
Using $\tilde{\phi}_{T,n}$, we obtain a quantile estimator as described in Section~\ref{sec:genQuan}:
\begin{enumerate}
\item The second derivative of the characteristic exponent can be estimated
by differentiating \eqref{eq:phiTilde} twice. We define
\begin{align}
\tilde\psi_n'(u) & :=-\frac{(2u+i)\F\tilde{O}(u+i)+u(iu-1)\F[x\tilde{O}](u+i)}{T(1-u(u-i)\F\tilde{O}(u+i))},\notag\\
\tilde{\psi}_{n}''(u) 
 & :=-\frac{(2\F\tilde{O}(u+i)+(4iu-2)\F[x\tilde{O}](u+i)-(u^{2}+iu)\F[x^{2}\tilde{O}](u+i))}{T(1-u(u-i)\F\tilde{O}(u+i))}\notag\\
 & \qquad-\frac{1}{T}\Big(\frac{(2u+i)\F\tilde{O}(u+i)+u(iu-1)\F[x\tilde{O}](u+i)}{1-u(u-i)\F\tilde{O}(u+i)}\Big)^{2}.\label{eq:psiTildePP} 
\end{align}
Using a kernel $K$ with bandwidth $h>0$ satisfying \eqref{eq:propKernel}, we obtain the density estimator 
\[
\tilde{\nu}_{h}(t):=-t^{-2}\F^{-1}\Big[\tilde{\psi}_{n}''(u)\F K(hu)\Big](t),\qquad t\neq0.
\] 

\item Integrating $\tilde\nu_h$, the estimator of the generalized distribution function
is given by 
\[
\tilde{N}_{h}(t)=-\int g_{t}(x)\F^{-1}\Big[\tilde{\psi}_{n}''(u)\F K(hu)\Big](x)\d x,\quad t\neq0,
\]
with $g_{t}$ from (\ref{eq:DefDistEst}). 

\item For $\tau>0$ the quantile estimators are defined as the minimum contrast estimators
\[
\tilde{q}_{\tau,h}^{\pm}:=\argmin_{t\in[\eta_{n},\infty)}|\tilde{N}_{h}(\pm t)-\tau|
\]
with threshold value $\eta_{n}\downarrow0$.
\end{enumerate}
These estimators are well defined on the event $A=\{\forall u\in[-1/h,1/h]:\tilde\phi_{T,n}(u)\neq0\}$ whose probability increases if $\tilde\phi_{T,n}$ concentrates around the true $\phi_T$. In an idealized model \cite{soehl:2010} shows $P(A)=1$. Note that $\tilde \nu_h$ is different from the finite activity estimator by \citet{belomestnyReiss2006}. The k-function estimator from \citet{trabs:2011} in the self-decomposable model relies on a similar idea but using only the first derivative of $\psi$.

Since we want to concentrate on the main aspects in the error analysis and to avoid technicalities, we will work in the idealized Gaussian white noise model which was considered by \citet{soehl2014}
as well. Assume that the noise levels of the observations \eqref{eq:Observations} are given by the values $\delta_j=\delta(x_j)$, $j=1,\dots,n$, of some function $\delta:\mathbb{R}\to\mathbb{R}_{+}$. The observed strike prices are assumed to be the quantiles $x_j=F^{-1}(j/(n+1))$, $j=1,\dots,n$, of a distribution with distribution function $F:\mathbb{R}\to[0,1]$ and density $f>0$. Incorporating the observation errors as well as their distribution, we define the general noise level
\begin{equation}\label{varrho}
  \varrho(x)=\delta(x)/\sqrt{f(x)}.
\end{equation}
For standard normal $(\eps_j)$ \cite{brownLow1996} have shown asymptotic equivalence in the sense of Le Cam of the nonparametric regression model~\eqref{eq:Observations} and the Gaussian white noise model 
$$\d Z(x)=O(x)\d x+n^{-1/2}\varrho(x)\d W(x)$$ 
with a two-sided Brownian motion $W$ for $x$ on a possibly growing bounded interval. The equivalence to more general error distributions follows from \cite{gramaNussbaum:2002}. More details on this equivalence can be found in \cite{soehl2014} and \citet[Supplement]{trabs:2011}. 

$Z$ is an empirical version of the antiderivative of $O$. In that sense we define $\F\tilde{O}(u):=\F[\d Z](u)=\F O(u)+n^{-1/2}\int e^{iux}\varrho(x)\d W(x)$ and analogously for $\F[x\tilde O]$ and $\F[x^2\tilde O]$. Owing to \eqref{eq:phiTilde}, the estimation
error of the $\tilde{\phi}_{T,n}$ is given by the Gaussian process
\[
\Phi_{n}(u):=(\tilde{\phi}_{T,n}-\phi_{T})(u)=u(u+i)\F[O-\tilde{O}](u+i)=n^{-1/2}u(u-i)\int e^{iux-x}\rho(x)\d W(x).
\]
While in the previous section the concentration of the estimator $\phi_{\Delta,n}$ around $\phi_\Delta$ was obtained by the i.i.d. structure of the increments, we use here the concentration of Gaussian measures. Applying Dudley's entropy theorem, we obtain the following path property of
$\Phi_n$. This lemma is in line with the results by \citet{soehl:2010}
and Proposition~1 in \citet{soehl2014}.
\begin{lemma}
\label{lem:GaussProc}Grant $\int(1+|x|)^{m}e^{-2x}\rho^{2}(x)\d x<\infty$
for some $m>4$. Then $\Phi_n$ is twice $L^{2}(P)$-differentiable
with derivatives
\begin{align}
\Phi_{n}^{(1)}(u):= & n^{-1/2}\int e^{iux-x}\big(2u-i+(iu^{2}+u)x\big)\rho(x)\d W(x),\label{eq:GaussProc}\\
\Phi_{n}^{(2)}(u):= & n^{-1/2}\int e^{iux-x}\big(2+2(2iu+1)x+u(i-u)x^{2}\big)\rho(x)\d W(x).\nonumber 
\end{align}
Moreover, $\Phi_{n}^{(0)}:=\Phi_{n},\Phi_{n}^{(1)}$ and $\Phi_{n}^{(2)}$ have versions
that are almost surely continuous and satisfy for any $U>0$
\[
\E\big[\|\Phi_{n}^{(k)}\|_{L^{\infty}[-U,U]}\big]=\mathcal O\big( n^{-1/2}U^{2}\sqrt{\log U}\big)\quad\text{for}\quad k=0,1,2.
\]

\end{lemma}
In the following we may use these almost surely continuous and bounded versions.
Let us first state a result on the uniform loss for the density estimator
$\tilde{\nu}_{h}$.
\begin{prop}
  \label{prop:UniformDensRates} Let $\alpha,\beta,s,r,R>0$, suppose the
  kernel satisfies (\ref{eq:propKernel}) with order $p\ge s$ and let
  $U\subset\R$ be a bounded, open set which is bounded away from zero.
  Suppose $\|(1\vee x^{2})e^{-x}\rho(x)\|_{\infty}<\infty$. Then we
  have for $n\to\infty$
  \begin{enumerate}
    \item with $h=h_{n}=\big(\frac{\log n}{n}\big)^{1/(2s+2T\alpha+5)}$ uniformly in $(\sigma^{2},\gamma,\nu)\in\mathcal{D}^{s}(\alpha,2,U,R)$ satisfying \eqref{eqMartCond} 
    \begin{align*}
    \sup_{t\in U}|\tilde{\nu}_{h}(t)-\nu(t)| & =\mathcal{O}_{P,\mathcal D^s}\Big(\Big(\frac{\log n}{n}\Big)^{s/(2s+2T\alpha+5)}\Big),
    \end{align*}
    \item with $h=h_{n}=(\frac{1}{4r}\log n)^{-1/\beta}$ uniformly in $(\sigma^{2},\gamma,\nu)\in\mathcal{E}^{s}(\beta,2,U,r,R)$ satisfying \eqref{eqMartCond}
    \[
    \sup_{t\in U}|\tilde{\nu}_{h}(t)-\nu(t)|=\mathcal{O}_{P,\mathcal E^s}\big((\log n)^{-s/\beta}\big).
    \]
  \end{enumerate}
\end{prop}
The convergence rates for the pointwise loss are the same without the logarithmic factor in (i). They coincide with the rates by \cite{belomestnyReiss2006} who have considered only the extreme cases: If $\sigma^2=0$ and $\nu$ is a finite measure, the pointwise risk converges with rate $n^{-s(2s+5)}$, and if $\sigma^2>0$, we obtain the rate $(\log n)^{-s/2}$. The estimator for the k-function by \cite{trabs:2011} achieves the same rate as the corresponding pointwise result in Proposition~\ref{prop:UniformDensRates}(i). Since in the two afore mentioned papers lower bounds have been proved and the logarithm is unavoidable for uniform loss, the above rates appear to be minimax optimal.

Recalling the function classes from \eqref{eq:Dtilde}, we obtain the following convergence rates for the quantile estimators $\tilde{q}_{\tau,h}^{\pm}.$
\begin{theorem}
  \label{thm:quantileFin}Let $\tau>0$ and $\alpha,\beta,s,\zeta,r,R>0,s'\in(-1,0]$.
  Suppose $\|(1\vee x^{2})e^{-x}\rho\|_{\infty}\lesssim1$ and $\int(1+|x|)^{m}e^{-2x}\rho^{2}(x)\d x<\infty$  for some $m>4$. Suppose the kernel satisfies (\ref{eq:propKernel})
  with order $p\ge s+1$ and let $\eta_{n}>0$ with $\eta_{n}^{-1}\lesssim\log n$.
  Then we obtain for $n\to\infty$:
  \begin{enumerate}
    \item with $h=h_{n}=n^{-1/(2s+2T\alpha+5)}$ uniformly in $\mathcal{\tilde{D}}_{\tau}^{s,s'}(\alpha,2,\zeta,\eta_n,R)$  satisfying \eqref{eqMartCond}  
    \begin{align*}
    |\tilde{q}_{\tau,h}^{\pm}-q_{\tau}^{\pm}| & =\mathcal{O}_{P,\mathcal{\tilde{D}}_{\tau}^{s,s'}}(n^{-(s+1)/(2s+2T\alpha+5)}),
    \end{align*}
    \item with $h=h_{n}=(\frac{1}{4r}\log n)^{-1/\beta}$ uniformly in $\mathcal{\tilde{E}}_{\tau}^{s,s'}(\beta,2,\zeta,\eta_n,r,R)$  satisfying \eqref{eqMartCond} 
    \[
    |\tilde{q}_{\tau,h}^{\pm}-q_{\tau}^{\pm}|=\mathcal{O}_{P, \mathcal{\tilde{E}}_{\tau}^{s,s'}}\big((\log n)^{-(s+1)/\beta}\big).
    \]
  \end{enumerate}
\end{theorem}
Compared to Theorem~\ref{thm:rateQuantiles}, the rates in (i) are always slower. In particular, the parametric rate can never be achieved. Heuristically, this is because we estimate a derived parameter of the state price density, which is basically the second derivative of the observed option function $O$. In the Fourier domain we see in the pricing formula~\eqref{pricingFormula} that $\F O$ decays two polynomial degrees faster than $\phi_T$ such that the ill-posedness of the statistical problem is larger. In the severely ill-posed case (ii) the rate is the same in both observation schemes since the rates are only logarithmically slow. The moment assumption is weakened to second moments, which are necessary for the identification identity~\eqref{eq:psiPP}. Instead the existence of fourth moments is implicitly imposed on the error distribution in the regression scheme. Although the theorem is stated for $\mathcal{\tilde{D}}_{\tau}^{s,s'}(\alpha,2,\zeta,\eta,R)$, the bounded variation of $x\nu$ is not needed here and could be dropped. 

\section{Data-driven choice of the bandwidth}\label{sec:lepski}

Of course the optimal bandwidth is not known to the practitioner.
To provide an adaptive method, we apply the approach by \cite{Lepski1990} by considering
a family of quantile estimators $\{\tilde{q}_{\tau,h}:h\in\mathcal{B}_{n}\}$
for an appropriate set of bandwidths $\mathcal{B}_{n}$. Following the
construction in \cite{dattnerEtAl2013}, we define for a constant
$L>1$ and a sequence $(N_{n})\subset\mathbb{N}$ satisfying $n^{-1}L^{N_{n}}\sim(\log n)^{-5}$
\[
h_{n,j}:=n^{-1}L^{j}\quad\text{ for }\quad j=0,\dots,N_{n}.
\]
To ensure that the density estimator $\tilde{\nu}_{h}$ is consistent
for any $h\in\mathcal{B}_{n}$, we choose the minimal bandwidth via
\begin{align*}
\tilde{j}_{n}:=&\min\Big\{ j=0,\dots,N_{n}:
&\frac{1}{2}\le \frac{(\log n)^{2}}{n^{1/2}}\Big(\int_{-1/h_{n,j}}^{1/h_{n,j}}\frac{(1+u^{4})\1_{\{|\tilde{\phi}_{T,n}(u)|\ge(1+|u|)^2/n^{1/2}\}}}{|\tilde{\phi}_{T,n}(u)|^{2}}\d u\Big)^{1/2}\le1\Big\}
\end{align*}
and define
\[
  \mathcal{B}_{n}:=\{h_{n,\tilde{j}_{n}},\dots,h_{n,N_{n}}\}.
\]
To choose the bandwidth from $\mathcal B_n$ which mimics the oracle bandwidth in Theorem~\ref{thm:quantileFin}, we have to estimate the standard deviation of the stochastic error. This problem is similar to the one considered by \cite{soehl2014} who has determined the asymptotic distribution of the finite activity estimators by \cite{belomestnyReiss2006} and has derived confidence sets. In fact, we will only estimate an upper bound of the standard deviation. This bound should be as sharp as possible to allow for a good finite sample behavior. The main problem is that any upper bound depends on unknown quantities which have to be estimated as well with a sufficiently fast convergence rate. Our bound will not depend on any asymptotic or function class specific constants. The stochastic error is dominated by its linearization and thus we define, cf. Lemma~\ref{lem:LinFin},
\begin{align}
  \tilde\Sigma_{n,h}^\pm:=&\frac{1}{2\pi n^{1/2}T}\Big(\|x^{2}e^{-x}\rho(x)\|_{\infty}\|\tilde\chi_{\tilde q_{\tau,h}^\pm}^{(2)}\|_{L^2}
  +\|xe^{-x}\rho(x)\|_{\infty}\|\tilde\chi_{\tilde q_{\tau,h}^\pm}^{(1)}\|_{L^2}+\|e^{-x}\rho(x)\|_{\infty}\|\tilde\chi_{\tilde q_{\tau,h}^\pm}^{(0)}\|_{L^2}\Big)\label{eq:SigmaTilde}
\end{align}
with auxiliary functions, for $u\in\R,t\neq0$,
\begin{align*}
  \tilde\chi_t^{(0)}(u)&:=\F g_{t}(-u)\F K(hu)\Big(u(u-i)\frac{T^2\tilde\psi_n'(u)^2-T\tilde\psi_n''(u)}{\tilde\phi_{T,n}(u)}+2T(i-2u)\frac{\tilde\psi_n'(u)}{\tilde\phi_{T,n}(u)}+2\tilde\phi_{T,n}^{-1}(u)\Big),\\
  \tilde\chi_t^{(1)}(u)&:=\F g_{t}(-u)\F K(hu)\Big((4iu+2)\tilde\phi_{T,n}^{-1}(u)-2Tu(iu+1)\frac{\tilde\psi_n'(u)}{\tilde\phi_{T,n}(u)}\Big),\\
  \tilde\chi_t^{(2)}(u)&:=u(i-u)\F g_{t}(-u)\F K(hu)\tilde\phi_{T,n}^{-1}(u).
\end{align*}
Note that $\tilde{\Sigma}_{n,h}^\pm$ are monotone decreasing in $h$. The magnitude of the stochastic error of $\tilde{q}_{\tau,h}^{\pm}$
can then be estimated by 
\begin{equation}\label{eq:V}
\tilde{V}_{n}^{\pm}(h):=\frac{(1+\delta)\sqrt{2\log\log n}\tilde{\Sigma}^\pm_{n,h}}{|\tilde{\nu}_{h}(\tilde{q}_{\tau,h}^{\pm})|}
\end{equation}
for any small $\delta>0$. Defining 
\[
\mathcal{U}_{h}^{\pm}:=[\tilde{q}_{\tau,h}^{\pm}-\tilde{V}_{n}^{\pm}(h),\tilde{q}_{\tau,h}^{\pm}+\tilde{V}_{n}^{\pm}(h)],
\]
the adaptive estimator is defined as 
\[
\tilde{q}_{\tau}^{\pm}:=\tilde{q}_{\tau,\tilde{h}^{\pm}}^{\pm}\quad\text{with}\quad\tilde{h}^{\pm}:=\max\{h\in\mathcal{B}_{n}:\bigcap_{\mu\le h,\mu\in\mathcal{B}_{n}}\mathcal{U}_{h}^{\pm}\neq\emptyset\}.
\]

\begin{theorem}
  \label{thm:quantileFinAdapt}Let $\tau>0$ and $\alpha,\beta,s,\zeta,r,R>0,s'\in(-1,0]$.
  Suppose $\|(1\vee x^{2})e^{-x}\rho\|_{\infty}\lesssim1$ and $\int(1+|x|)^{m}e^{-2x}\rho^{2}(x)\d x<\infty$ for some $m>4$. Suppose the kernel satisfies (\ref{eq:propKernel})
  with order $p\ge s+1$ and let $\eta_{n}>0$ with $\eta_{n}^{-1}\lesssim\log n$.
  Then we obtain for $n\to\infty$: 
  \begin{enumerate}
    \item uniformly in $(\sigma^2,\gamma,\nu)\in\mathcal{\tilde{D}}_{\tau}^{s,s'}(\alpha,2,\zeta,\eta_n,R)$ satisfying \eqref{eqMartCond}
    \begin{align*}
    |\tilde{q}_{\tau}^{\pm}-q_{\tau}^{\pm}| & =\mathcal{O}_{P,\mathcal{\tilde{D}}_{\tau}^{s,s'}}\big(((\log\log n)n^{-1})^{(s+1)/(2s+2T\alpha+5)}\big),
    \end{align*}

    \item uniformly in $(\sigma^2,\gamma,\nu)\in\mathcal{\tilde{E}}_{\tau}^{s,s'}(\beta,2,\zeta,\eta_n,r,R)$ satisfying \eqref{eqMartCond}
    \[
    |\tilde{q}_{\tau}^{\pm}-q_{\tau}^{\pm}|=\mathcal{O}_{P,\mathcal{\tilde{E}}_{\tau}^{s,s'}}\big((\log n)^{-(s+1)/\beta}\big).
    \]
  \end{enumerate}
\end{theorem}
In the mildly ill-posed case (i) the adaptive method looses a $\log\log n$-factor compared to the oracle choice in Theorem~\ref{thm:quantileFin}. Note that our loss function is bounded thus we loose only a $\log\log n$-factor instead of the $\log n$-factor, which would appear, for instance, for the mean squared error. In view of \cite{spokoiny1996} this payment for adaptivity is unavoidable. In the severely ill-posed case (ii) the rates are already logarithmically slow such that the adaptive method causes no additional loss.

\section{Simulations and real data }\label{sec:sim}

We will illustrate the quantile estimation method in simulations from the CMGY model introduced by \cite{Carr:2002}. The driving L\'evy process of the asset  is tempered stable and may have a diffusion component. For parameters $C>0,M,G\ge0$ and $Y<2$ the L\'evy measure in the CMGY model is given by the Lebesgue density 
\[
  \nu_{cgmy}(x)=\begin{cases}
           C|x|^{-1-Y}e^{-G|x|},\quad &x<0,\\
           Cx^{-1-Y}e^{-Mx},\quad &x>0.
         \end{cases}
\]
With a fifth parameter $\sigma\ge0$ the characteristic triplet of the underlying L\'evy process is $(\sigma^2,\gamma,\nu_{cgmy})$ where the drift is determined by the martingale condition \eqref{eqMartCond}. For the simulations we set $C=1, G=5, M=8, Y=0.5$ and $\sigma=0.1$ which appears to be a realistic choice in view of the empirical results by \cite{Carr:2002}. The riskless interest rate is chosen as $r=0.06$. The design points $x_j,j=1,\dots,n,$ are constructed as $j/(n+1)$-quantiles of a $\mathcal N(0,1/2)$ distribution. We simulate $n=100$ option prices with time to maturity $T=0.25$ corresponding to three months. According to a rule of thumb by \citet[cf.][p. 439]{contTankov:2004b}, the local noise levels $\delta_j$ are chosen as 1\% of the observed prices $O(x_j)$ which is later assumed for the real data, too.

\begin{table}\centering
  \begin{tabular}{ccccccc}
    \hline\hline
     \multicolumn{3}{c}{$10^2\cdot$RMSE}&\multicolumn{2}{c}{$\tilde q_\tau^-$} & \multicolumn{2}{c}{$\tilde q_\tau^+$}\\
    $\tau$ & $q_\tau^-$ & $q_\tau^+$ & oracle & adaptive  & oracle & adaptive \\\hline
    0.5 & 0.1778 & 0.1241 & 0.346 & 4.806 & 0.444 & 2.246\\
    1.0 & 0.1201 & 0.0868 & 0.297 & 1.396 & 0.361 & 0.741\\
    1.5 & 0.0929 & 0.0665 & 0.185 & 0.890 & 0.434 & 0.869\\
    2.0 & 0.0726 & 0.0563 & 0.275 & 0.867 & 0.314 & 0.670\\
    2.5 & 0.0624 & 0.0461 & 0.233 & 0.652 & 0.424 & 0.694\\\hline\hline
   \end{tabular}
  \caption{Empirical RMSE (multiplied by 100) of the quantile estimators $\tilde q_\tau^\pm$ from 1000 Monte Carlo simulations of the CGMY model.}\label{tab:cgmy}
\end{table}

In order to apply the estimation procedure, we have to choose some parameters. The truncation value is set to $\eta=0.02$. To construct the bandwidth set $\mathcal B_n$, we take $L=1.1$. To compute $\tilde \Sigma_{n,h}^{\pm}$, we need the noise function $\varrho$ from \eqref{varrho} which depends on the unknown density $f$ of the distribution of the strikes. It can be estimated from the observation points $(x_j)_{j=1,\dots,n}$ using some standard density estimation method. As in \cite{soehlTrabs2014} we will apply a triangular kernel estimator, where the bandwidth is chosen by Silverman's rule of thumb.

To assess the performance of the estimation procedure, we compare the adaptive choice of the bandwidth to the oracle bandwidth, meaning that $h$ is chosen such that the empirical root mean squared error (RMSE) is minimized. Our simulation results are summarized in Table~\ref{tab:cgmy} for $\tau\in\{0.5,1.0,1.5,2.0,2.5\}$. Although the RMSE of the adaptive method is larger than the oracle choice, the method achieves reasonable estimation errors. Note that the sample size is relatively small. For $\tau=0.5$ the RMSE is about 2.7\% of the true quantile for $\tilde q_{0.5}^-$ and 1.8\% for $\tilde q_{0.5}^+$. For larger values of $\tau$ the RMSE decays to an order of 1\% of the true quantiles. The reason is that small values of $\tau$ correspond to rare large jumps such that the jump density is small. Consequently, the estimation error \eqref{eq:errorRep} is large. Since the stochastic estimation error has to be estimated by $\tilde V_n^\pm$ from \eqref{eq:V}, this effect is more severe for the Lepski method. 

\begin{figure}
  \includegraphics[width=14cm]{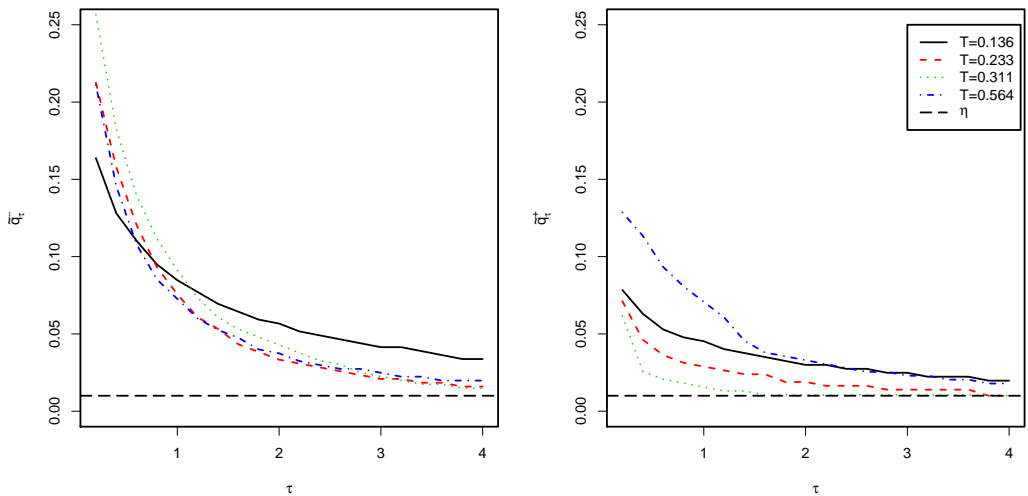}
  \caption{Estimated generalized quantiles of negative jumps \textit{(left)} and positive jumps \textit{(right)} based on option prices from May 29, 2008, with four maturities $T$.}\label{fig:quantilesDax}
\end{figure}
Let us finally apply the estimation method to prices of DAX options from May 29, 2008. This data set\footnote{provided by the SFB 649 ``Economic Risk''} has already been studied by \cite{soehlTrabs2014}. Figure~\ref{fig:quantilesDax} shows the estimated quantiles $\tilde q_\tau^\pm$ for four different maturities between two and seven months and for $\tau\in\{0.2,0.4,\dots,4\}$. Due to this finer grid, the threshold value is set to $\eta=0.01$. As postulated by the stylized facts on financial data, negative jumps have a higher activity. Roughly, the intensity of small jumps is larger for short maturities while the tails are more heavy for longer maturities.

\section{Proofs}\label{sec:LevyProofs}

\subsection{Error analysis for Section~\ref{sec:quantileEst}}

To simplify the notation we will frequently use the definition $I_{h}:=[-\frac{1}{h},\frac{1}{h}]$. Note that the indicator function $\1_{\{|\phi_{\Delta,n}(u)|\ge(\Delta n)^{-1/2}\}}$ in the definition of $\hat\psi_n''$ equals one on $I_h$ with probability converging to one. This follows exactly from Lemma~5.1 by \cite{dattnerEtAl2013} and the bandwidth choices which we will consider.

\subsubsection{Drift and remainder }

All estimators are constructed based on the estimator $\hat\psi_n''(u)=\Delta^{-1}(\log(\phi_{\Delta,n}(u))''$ of $\psi''$ from \eqref{eq:psiPP}. In particular, $\hat\nu_h,\hat N_h$ and $\hat q_{\tau,h}$ only depend on the observations via $\hat\psi_n''$. As shown by \citet[Lem. 10]{nicklEtAl2013} the drift has no effect on the estimators: 
\begin{lemma}
\label{lem:drift}Let $X_{k}:=Y_{k}-\Delta\gamma$ for $k=1,\dots,n.$
Then $X_{k}$ are distributed according to an infinitely divisible
distribution with characteristic triplet $(\sigma^{2},0,\nu)$. Denoting
the estimator $\hat\psi_n''$ based on $(X_{k})$ and $(Y_{k})$
by $\hat\psi''_{X,n}$ and $\hat{\psi}''_{Y,n}$, respectively, it holds $\hat\psi''_{X,n}=\hat\psi''_{Y,n}$
for all $u\in\R.$
\end{lemma}

Using 
\begin{equation}
(\phi_{\Delta}^{-1})'=-\Delta\psi'\phi_{\Delta}^{-1},\qquad(\phi_{\Delta}^{-1})''=-\Delta\big(\psi''-\Delta(\psi'){}^{2}\big)\phi_{\Delta}^{-1},\label{eq:PhiDerivatives}
\end{equation}
the estimation error $\psi''-\hat{\psi}''_{n}=-\Delta^{-1}\log(\phi_{\Delta,n}/\phi_{\Delta})''$ can be linearized similarly to Proposition~21 in \citet{nicklEtAl2013}. 
\begin{lemma}
\label{lem:Remainder}Let $\int x^{4+\delta}\nu(\d x)<\infty$ for
some $\delta>0$. For $h,\Delta\in(0,1)$ satisfying \\$n^{-1/2}(\log h^{-1})^{(1+\delta)/2}\|\phi_{\Delta}^{-1}\|_{L^{\infty}(I_h)}\to0$
as $n\to\infty$, it holds
\begin{align*}
 & \sup_{|u|\le h^{-1}}\left|\hat{\psi}_{n}''(u)-\psi''(u)-\Delta^{-1}(\phi_{\Delta}^{-1}(\phi_{\Delta,n}-\phi_{\Delta}))''(u)\right|\\
= & \mathcal O_{P}\big((\Delta\|\psi'\big\|_{L^{\infty}(I_h)}+\Delta^{3/2}\|\psi'\big\|_{L^{\infty}(I_h)}^{2}+1)n^{-1}\Delta^{-1/2}\log(h^{-1})^{1+\delta}\|\phi_{\Delta}^{-1}\|_{L^{\infty}(I_h)}^{2}\big).
\end{align*}
\end{lemma}
\begin{proof}
Setting $F(y)=\log(1+y)$ and $\eta=(\phi_{\Delta,n}-\phi_{\Delta})/\phi_{\Delta}$,
we use $(F\circ\eta)''(u)=F'(\eta(u))\eta''(u)+F''(\eta(u))\eta'(u)^{2}$
to obtain 
\begin{align}
|(F\circ\eta)''(u)-\eta''(u)| & \lesssim\|F''\|_{\infty}\Big(\eta(u)\eta''(u)+\eta'(u)^{2}\Big).\label{eq:MarkusTrick}
\end{align}
On the event $\Omega_{n}:=\{\sup_{|u|\le1/h}|(\phi_{\Delta,n}-\phi_{\Delta})(u)/\phi_{\Delta}(u)|\le1/2\}$
we thus obtain 
\begin{align*}
 & \sup_{\left|u\right|\le h^{-1}}\big|\log(\phi_{\Delta,n}/\phi_{\Delta})''(u)-(\phi_{\Delta}^{-1}(\phi_{\Delta,n}-\phi_{\Delta}))''(u)\big|\\
 & \qquad\lesssim\|\eta\|_{L^{\infty}(I_h)}\|\eta''\|_{L^{\infty}(I_h)}+\|\eta'\|_{L^{\infty}(I_h)}^{2}.
\end{align*}
To estimate $\|\eta^{(k)}\|_{L^{\infty}(I_h)},k=0,1,2$, we use \eqref{eq:PhiDerivatives} and
$\left|\psi''(u)\right|\lesssim1$ to obtain 
\begin{align*}
\sup_{u\in I_h}\left|(\phi_{\Delta}^{-1})'(u)\right| & \lesssim\Delta\|\psi'\big\|_{L^{\infty}(I_h)}\|\phi_{\Delta}^{-1}\|_{L^{\infty}(I_h)},\\
\sup_{u\in I_h}\left|(\phi_{\Delta}^{-1})''(u)\right| & \lesssim(\Delta^{2}\|\psi'\big\|_{L^{\infty}(I_h)}^{2}+\Delta)\|\phi_{\Delta}^{-1}\|_{L^{\infty}(I_h)}.
\end{align*}
Applying Theorem 1 by \citet{kappusReiss2010} and the moment assumption
on $\nu$, we obtain for $k=0,1,2$ and $\delta>0$ 
\[
\big\|(\phi_{\Delta,n}-\phi_{\Delta})^{(k)}\big\|_{L^{\infty}(I_h)}=\mathcal O_{P}(n^{-1/2}\Delta^{(k\wedge1)/2}(\log h^{-1})^{(1+\delta)/2}).
\]
This yields 
\begin{align}
\|\eta\|_{L^{\infty}(I_h)} & =\mathcal O_{P}\big(n^{-1/2}(\log h^{-1})^{(1+\delta)/2}\|\phi_{\Delta}^{-1}\|_{L^{\infty}(I_h)}\big),\label{eq:supEta}\\
\|\eta'\|_{L^{\infty}(I_h)} & =\mathcal O_{P}\big((\Delta\|\psi'\big\|_{L^{\infty}(I_h)}+\Delta^{1/2})n^{-1/2}(\log h^{-1})^{(1+\delta)/2}\|\phi_{\Delta}^{-1}\|_{L^{\infty}(I_h)}\big),\nonumber \\
\|\eta''\|_{L^{\infty}(I_h)} & =\mathcal O_{P}\big((\Delta^{3/2}\|\psi'\big\|_{L^{\infty}(I_h)}+\Delta^{2}\|\psi'\big\|_{L^{\infty}(I_h)}^{2}+\Delta^{1/2})\nonumber \\
 & \hspace{3em}\times n^{-1/2}(\log h^{-1})^{(1+\delta)/2}\|\phi_{\Delta}^{-1}\|_{L^{\infty}(I_h)}\big).\nonumber 
\end{align}
The bound (\ref{eq:supEta}) and $n^{-1/2}(\log h^{-1})^{(1+\delta)/2}\|\phi_{\Delta}^{-1}\|_{L^{\infty}(I_h)}\to0$,
 yield $P(\Omega_{n})\to1$ which implies
the assertion. 
\end{proof}

\subsubsection{Convergence rates for distribution function estimation}
Because it is a bit more difficult to derive the convergence rates for the distribution function estimator than for the density estimator, we study this problem first. The corresponding results for $\hat\nu_n$
can the be proved analogously. 

We decompose the estimation error of the distribution function estimator
$\hat{N}_{h}$ into
\begin{align}
\hat{N}_{h}(t)-N(t)= & \int g_{t}(x)\big(K_{h}\ast\big(y^{2}\nu\big)-x^{2}\nu\big)(\d x)\nonumber \\
 & +\int g_{t}(x)\F^{-1}\Big[\Big(\psi''(u)-\hat{\psi}''_{n}(u)\Big)\F K(hu)\Big](x)\d x+\sigma^{2}\int g_{t}(x)K_{h}(x)\d x\nonumber \\
=: & B_{n}(t)+S_{n}(t)+V_{n}(t),\label{eq:Decomp}
\end{align}
where $B_{n}$ is the deterministic error term, $S_{n}$ is the stochastic
error term and $V_{n}$ is the error due to the unknown volatility
$\sigma^{2}$. The error term $V_{n}$ is negligible:
\begin{lemma}
\label{lem:vola}Grant Assumption (\ref{eq:propKernel}) on the kernel
with order $p\ge s+1$. Then $V_{n}(t),t\neq0,$ as defined in (\ref{eq:Decomp})
satisfies $|V_{n}(t)|\lesssim\sigma^{2}|t|^{-s-3}h^{s+1}.$\end{lemma}
\begin{proof}
We estimate 
\begin{align*}
|V_{n}|\le\sigma^{2}\int|g_{t}(x)K_{h}(x)|\d x\le & \sigma^{2}\|x^{-s-1}g_{t}(x)\|_{\infty}\|x^{s+1}K_{h}(x)\|_{L^{1}}\\
= & \sigma^{2}|t|^{-s-3}h^{s+1}\|x^{s+1}K(x)\|_{L^{1}}.\qedhere
\end{align*}

\end{proof}
For the bias we apply the following:
\begin{prop}\label{prop:bias}
  Suppose $\|x^{2}\nu\|_{L^{1}}<\infty$ and let $U\subset\R$
  be an open set. If $\nu$ admits a Lebesgue
  density on $U$ in $C^{s}(U)$ for some $s>-1$ and if the kernel satisfies
  (\ref{eq:propKernel}) with order $p\ge s+1$, then 
  \[
  |B_{n}(t)|=\Big|\int g_{t}(x)\big(K_{h}\ast\big(y^{2}\nu(\d y)\big)-x^{2}\nu\big)(\d x)\Big|\lesssim(|t|^{-s-4}\vee1)h^{s+1},\quad\text{for all}\quad t\in U.
  \]
\end{prop}
\begin{proof}
Without loss of generality let $t<0.$ Using Fubini's theorem, we
rewrite 
\begin{align}
\int_{\R}g_{t}(x)\big(K_{h}\ast\big(y^{2}\nu(\d y)\big)-x^{2}\nu\big)(\d x)=K_{h}\ast g_{t}(-\bull)\ast(x^{2}\nu)(0)-g_{t}(-\bull)\ast(x^{2}\nu)(0).\label{eqRewrBias}
\end{align}
Denoting $\overline{N}(t)=\int_{-\infty}^{t}\nu(\d x),t<0,$ integration
by parts yields 
\begin{align*}
g_{t}(-\bull)\ast(x^{2}\nu)(y) & =\int_{-\infty}^{t+y}\frac{x^{2}}{(x-y)^{2}}\nu(\d x)\\
 & =\frac{(t+y)^2}{t^2}\overline{N}(t+y)+\int_{-\infty}^{0}\frac{2(x+t+y)y}{(x+t)^{3}}\overline{N}(x+t+y)\d x\\
 & =\frac{(t+y)^2}{t^2}\overline{N}(t+y)+2y\left(((\bull+t)^{-2}\1_{(-\infty,0]})\ast\overline{N}(t-\bull)\right)(-y) \\
  &\qquad+2y^2\left(((\bull+t)^{-3}\1_{(-\infty,0]})\ast\overline{N}(t-\bull)\right)(-y).
\end{align*}
From this representation we see for some sufficiently small $\delta\in(0,|t|)$
that $g_{t}(-\bull)\ast(x^{2}\nu)\in C^{s+1}((-\delta/2,\delta/2))$
owing to $\overline{N}(t+\bull)\in C^{s+1}((-\delta,\delta))$ and $(\bull+t)^{-3}\1_{(-\infty,0]}\in C^{s+1}((-\delta/2,\delta/2)^{c})$.
The corresponding H\"older norm of the latter is of order $t^{-s-4}$.
With a standard Taylor expansion argument and applying the order of
the kernel, the approximation error is of the order $(|t|^{-s-4}\vee1)h^{s+1}$ by
(\ref{eqRewrBias}). 
\end{proof}

Lemma~\ref{lem:Remainder} motivates the following definition of
the linearized stochastic error term 
\begin{equation}
L_{\Delta,n}(t):=-\Delta^{-1}\int g_{t}(x){\cal F}^{-1}[\F K(h\bull)\big(\phi_{\Delta}^{-1}(\phi_{\Delta,n}-\phi_{\Delta})\big)''](x)\d x.\label{eq:linError}
\end{equation}
Using \eqref{eq:PhiDerivatives} and defining the regularized Fourier multiplier
\[
m_{\Delta,h}:=\frac{\F K(h\bull)}{\phi_{\Delta}},
\]
we decompose the linearized stochastic error further into 
\begin{align}
L_{\Delta,n}(t)= & -\frac{1}{\Delta}\int g_{t}(x){\cal F}^{-1}[\F K(h\bull)\big(\phi_{\Delta}^{-1}(\phi_{\Delta,n}-\phi_{\Delta})''+2(\phi_{\Delta}^{-1})'(\phi_{\Delta,n}-\phi_{\Delta})'\nonumber \\
 & \quad+(\phi_{\Delta}^{-1})''(\phi_{\Delta,n}-\phi_{\Delta})\big)](x)\d x\nonumber \\
= & \underbrace{-\frac{1}{\Delta}\int g_{t}(x){\cal F}^{-1}[m_{\Delta,h}(\phi_{\Delta,n}''-\phi_{\Delta}'')](x)\d x}_{=:M_{\Delta,n}(t)}\notag\\
&\qquad+2\int g_{t}(x){\cal F}^{-1}[m_{\Delta,h}\psi'(\phi_{\Delta,n}'-\phi_{\Delta}')](x)\d x\nonumber \\
 & \qquad+\int g_{t}(x){\cal F}^{-1}[m_{\Delta,h}\big(\psi''-\Delta(\psi'){}^{2}\big)(\phi_{\Delta,n}-\phi_{\Delta})](x)\d x.\label{eq:DecompLin}
\end{align}
In the following, we will refer to $M_{\Delta,n}$ as the main stochastic
error term.
\begin{prop}\label{prop:variance}
  Let $U\subset\R$ be an open set, $t\in U$
  and $\alpha,\beta,s,r,R>0,m>4$. Let the kernel satisfy (\ref{eq:propKernel})
  with $p\ge1$. Then we obtain for $h,\Delta\in(0,1)$ and $\eta\in(0,1)$
  \begin{align*}
  \sup_{(\sigma^{2},\gamma,\nu)\in\mathcal{D}^{s}(\alpha,m,U,R)}\E\Big[\sup_{|t|>\eta}|L_{\Delta,n}-M_{\Delta,n}|(t)\Big] & \lesssim\eta^{-2}n{}^{-1/2}\|(1+|u|)^{\Delta\alpha-1}\|_{L^{2}(I_h)},\\
  \sup_{(\sigma^{2},\gamma,\nu)\in\mathcal{E}^{s}(\beta,m,U,r,R)}\E\Big[\sup_{|t|>\eta}|L_{\Delta,n}-M_{\Delta,n}|(t)\Big] & \lesssim\frac{\log(h^{-1})+\Delta^{1/2}h^{-1}+\Delta h^{-2}}{\eta^{2}n^{1/2}}e^{r\Delta h^{-\beta}}.
  \end{align*}
  Moreover,
  \[
  \sup_{(\sigma^{2},\gamma,\nu)\in\mathcal{E}^{s}(\beta,m,U,r,R)}\E\Big[\sup_{|t|>\eta}|M_{\Delta,n}(t)|\Big]\lesssim\eta^{-2}(n\Delta)^{-1/2}\log(h^{-1})e^{r\Delta h^{-\beta}}.
  \]
\end{prop}
\begin{remark}
Note that
\begin{align*}
\|(1+|u|)^{\Delta\alpha-1}\|_{L^{2}(I_h)} & \lesssim\begin{cases}
h^{-\Delta\alpha+1/2} & \text{for }\Delta\alpha>1/2,\\
(\log h^{-1})^{1/2} & \text{for }\Delta\alpha=1/2,\\
1 & \text{for }\Delta\alpha\in(0,1/2).
\end{cases}
\end{align*}
\end{remark}
\begin{proof}
Due to $\|g_{t}\|_{L^{1}}=|t|^{-1},\|g_{t}\|_{BV}\le2t^{-2}$, we
obtain $|\F g_{t}(u)|\lesssim(t^{-1}\vee t^{-2})(1+|u|)^{-1},u\in\R$.
Using Plancherel's identity, we estimate for $|t|>\eta$
\begin{align*}
|L_{\Delta,n}-M_{\Delta,n}|(t)\le & \pi^{-1}\Big|\int\F g_{t}(-u)m_{\Delta,h}(u)\psi'(u)(\phi_{\Delta,n}-\phi_{\Delta})'(u)\d u\Big|\\
 & +(2\pi)^{-1}\Big|\int\F g_{t}(-u)m_{\Delta,h}(u)\big(\psi''(u)-\Delta\psi'(u){}^{2}\big)(\phi_{\Delta,n}-\phi_{\Delta})(u)\d u\Big|\\
\lesssim & \eta^{-2}\int(1+|u|)^{-1}\big|m_{\Delta,h}(u)\psi'(u)(\phi_{\Delta,n}-\phi_{\Delta})'(u)\big|\d u\\
 & +\eta^{-2}\int(1+|u|)^{-1}\big|m_{\Delta,h}(u)\psi''(u)(\phi_{\Delta,n}-\phi_{\Delta})(u)\big|\d u\\
 & +\eta^{-2}\Delta\int(1+|u|)^{-1}\big|m_{\Delta,h}(u)\psi'(u){}^{2}(\phi_{\Delta,n}-\phi_{\Delta})(u)\big|\d u.
\end{align*}
Due to $\E[Y_{1}^{2l}]\lesssim\Delta^{l\wedge1},l=0,1,2$, we have, moreover,
\begin{equation}
\sup_{u\in\R}\E[(\phi_{\Delta,n}^{(l)}-\phi_{\Delta}^{(l)})^{2}(u)]\lesssim n^{-1}\Delta^{(l\wedge1)},\quad\text{for }l=0,1,2.\label{eq:phiMom}
\end{equation}
Fubini's theorem and Jensen's inequality yield
\begin{align}
 & \E\Big[\sup_{|t|>\eta}|L_{\Delta,n}-M_{\Delta,n}|(t)\Big]\nonumber \\
\lesssim & \eta^{-2}\Big(\int(1+|u|)^{-1}\big|m_{\Delta,h}(u)\psi'(u)\big|\E\big[(\phi'_{\Delta,n}-\phi'_{\Delta})^{2}(u)\big]^{1/2}\d u\nonumber \\
 & \qquad+\int(1+|u|)^{-1}\big|m_{\Delta,h}(u)\psi''(u)|\E[(\phi_{\Delta,n}-\phi_{\Delta})^{2}(u)]^{1/2}\d u\nonumber \\
 & \qquad+\Delta\int(1+|u|)^{-1}\big|m_{\Delta,h}(u)\psi'(u){}^{2}|\E[(\phi_{\Delta,n}-\phi_{\Delta})^{2}(u)]^{1/2}\d u\Big)\nonumber \\
\lesssim & n^{-1/2}\eta^{-2}\Big(\Delta^{1/2}\big\|(1+|u|)^{-1}m_{\Delta,h}(u)\psi'(u)\big\|_{L^{1}}+\big\|(1+|u|)^{-1}m_{\Delta,h}(u)\psi''(u)\big\|_{L^{1}}\nonumber \\
 & \qquad+\Delta\big\|(1+|u|)^{-1}m_{\Delta,h}(u)\psi'(u){}^{2}\big\|_{L^{1}}\Big).\label{eq:L-M}
\end{align}
To deal with $(\sigma^{2},\gamma,\nu)\in\mathcal{D}^{s}(\alpha,m,U,R)$,
we note that the assumptions $\|x^{4}\nu\|_{L^{1}}<\infty$ and $\|x\nu\|_{\infty}<\infty$
imply $x\nu,x^{2}\nu\in L^{1}(\R)\cap L^{2}(\R)$. Moreover, $\sigma^{2}=0$
and by Lemma~\ref{lem:drift} we can assume $\gamma_0=\gamma-\int x\d\nu=0$ such that  $i\psi'=-\F[x\nu]$ and $\psi''=-\F[x^{2}\nu]$. Therefore,
we estimate (\ref{eq:L-M}) with use of the Cauchy--Schwarz inequality
\begin{align}
 & \E\big[\sup_{|t|>\eta}|L_{\Delta,n}-M_{\Delta,n}|(t)\big]\nonumber \\
\lesssim & n^{-1/2}\eta^{-2}\Big(\Delta^{1/2}\big\|(1+|u|)^{-1}m_{\Delta,h}(u)\|_{L^2}\|\psi'\big\|_{L^{2}}+\big\|(1+|u|)^{-1}m_{\Delta,h}(u)\|_{L^{2}}\|\psi''\big\|_{L^{2}}\nonumber \\
 & \qquad+\Delta\big\|(1+|u|)^{-1}m_{\Delta,h}(u)\|_{L^{2}}\|\psi'\big\|_{L^{2}}\|\psi'\|_{\infty}\Big)\nonumber \\
\lesssim & n^{-1/2}\eta^{-2}\big(\Delta^{1/2}\|x\nu\|_{L^{2}}+\|x^{2}\nu\|_{L^{2}}+\Delta\|x\nu\|_{L^{2}}\|x\nu\|_{L^{1}}\big)\big\|(1+|u|)^{\Delta\alpha-1}\|_{L^{2}(I_h)},\label{eq:nonCritMild}
\end{align}
which yields the assertion for $\mathcal{D}^{s}(\alpha,m,U,R)$.

Now we consider the case $(\sigma^{2},\gamma,\nu)\in\mathcal{E}^{s}(\beta,m,U,r,R)$.
The exponential decay of $\phi_{\Delta}$, the properties of $K$,
$|\psi'(u)|\lesssim1+|u|$ and $|\psi''(u)|\lesssim1$ yield
\begin{align*}
\Delta^{1/2}\big\|(1+|u|)^{-1}m_{\Delta,h}(u)\psi'(u)\big\|_{L^{1}}\lesssim & \Delta^{1/2}\big\| m_{\Delta,h}\big\|_{L^{1}}\\
\lesssim&\Delta^{1/2}\int_{-1/h}^{1/h}e^{r\Delta|u|^{\beta}}\d u\lesssim\Delta^{1/2}h^{-1}\exp\big(r\Delta h^{-\beta}\big),\\
\big\|(1+|u|)^{-1}m_{\Delta,h}(u)\psi''(u)\|_{L^{1}}\lesssim & \big\|(1+u)^{-1}m_{\Delta,h}(u)\big\|_{L^{1}}\lesssim\log(h^{-1})\exp\big(r\Delta h^{-\beta}\big),\\
\Delta\big\|(1+|u|)^{-1}m_{\Delta,h}(u)\psi'(u){}^{2}\big]\big\|_{L^{1}}\lesssim & \Delta\big\|(1+u)m_{\Delta,h}(u)\big\|_{L^{1}}\lesssim\Delta h^{-2}\exp\big(r\Delta h^{-\beta}\big).
\end{align*}
We conclude the claimed estimate for $|L_{\Delta,n}-M_{\Delta,n}|$ in $\mathcal{E}^{s}(\beta,m,U,r,R)$
by plugging these estimates into (\ref{eq:L-M}). Similarly, we estimate
\begin{align*}
\E\Big[\sup_{|t|>\eta}|M_{\Delta,n}(t)|\Big] & \lesssim\eta^{-2}\Delta^{-1}\E\Big[\int(1+|u|)^{-1}m_{\Delta,h}(u)(\phi_{\Delta,n}''-\phi_{\Delta}'')(u)\d u\Big]\\
 & =\eta^{-2}(n\Delta)^{-1/2}\|(1+|u|)^{-1}m_{\Delta,h}(u)\|_{L^{1}}\\
 &\lesssim\eta^{-2}(n\Delta)^{-1/2}\log(h^{-1})\exp\big(r\Delta h^{-\beta}\big).\qedhere
\end{align*}

\end{proof}

For the main stochastic error term in the mildly ill-posed case we
will need the following concentration result:
\begin{prop}
\label{prop:mainStochErr}Let $U\subset\R$ be an open set, $t\in U$
and $\alpha,s,R>0,m>4$ and let the kernel satisfy (\ref{eq:propKernel})
with $p\ge1$. Then there is some $c>0$ such that for any $h,\Delta,\eta\in(0,1)$
and $\kappa_{0}>0$ 
\begin{align*}
\sup_{(\sigma^{2},\gamma,\nu)\in\mathcal{D}^{s}(\alpha,m,U,R)}\sup_{|t|>\eta} & P\Big(|M_{\Delta,n}(t)|>\kappa_{0}(\eta^{-1}\vee1)(n\Delta)^{-1/2}\|(1+|u|)^{\Delta\alpha-1}\|_{L^{2}(I_h)}\Big)\\
 & \hspace{6em}\le2\exp\Big(-\frac{c\kappa_{0}^{2}}{1+\kappa_{0}(\eta^{-1}\vee1)(hn\Delta)^{-1/2}}\Big).
\end{align*}
\end{prop}
\begin{proof}
We represent $M_{\Delta,n}(t)$ as sum of i.i.d. random variables
via 
\[
M_{\Delta,n}(t)=\sum_{k=1}^{n}(\xi_{k}(t)-\E[\xi_{k}(t)]),\quad\xi_{k}(t):=(n\Delta)^{-1}\int g_{t}(x)\F^{-1}\big[m_{\Delta,h}(u)Y_{k}^{2}e^{iuY_{k}}\big](x)\d x.
\]
Applying Plancherel's identity, $\xi_{k}(t)$ can be rewritten as
\begin{align*}
\xi_{k}(t) & =\frac{1}{2\pi n\Delta}Y_{k}^{2}\int\F g_{t}(-u)m_{\Delta,h}(u)e^{iuY_{k}}\d u=\frac{1}{n\Delta}Y_{k}^{2}\F^{-1}[\F g_{t}(-\bull)m_{\Delta,h}](-Y_{k}).
\end{align*}
To estimate $\Var(\xi_{k})\le\E[\xi_{k}^{2}]$, we use $g_{t}\in BV(\R)$
to decompose $g_{t}=g_{t}^{s}+g_{t}^{c}$ into a singular component
and a continuous component satisfying for $r>\Delta\alpha$
\begin{align*}
\max\left\{ \F g_{t}^{s}(u),\F[xg_{t}^{s}](u),\F[x^{2}g_{t}^{s}](u)\right\} &\lesssim(t^{-2}\vee1)(1+|u|)^{-1},\\
\max\Big\{\sup_{|t|>\eta}\|g_{t}^c\|_{C^r},\sup_{|t|>\eta}\|x^2g_{t}^c\|_{C^r}\Big\}&\lesssim1.
\end{align*}
This allows to decompose 
\begin{align}
\E[\xi_{k}^{2}] & \le\frac{2}{(n\Delta)^{2}}\Big(\E\Big[Y_{1}^{4}\F^{-1}\big[\F g_{t}^{s}(-\bull)m_{\Delta,h}\big](-Y_{1})^{2}\Big]\notag\\
&\qquad\qquad+\E\Big[Y_{1}^{4}\F^{-1}\big[\F g_{t}^{c}(-\bull)m_{\Delta,h}\big](-Y_{1})^{2}\Big]\Big) =:\frac{2}{(n\Delta)^{2}}(E_{s}+E_{c}).\label{eq:VarM}
\end{align}
To estimate $E_{c}$ in (\ref{eq:VarM}), we apply the Fourier multiplier
Theorem~5 in \cite{trabs2014} to see that 
\[
  \|\F^{-1}\big[\F g_{t}^{c}(-\bull)m_{\Delta,h}\big]\|_{\infty}\lesssim\|g_{t}^{c}\ast K_{h}\|_{C^{r}}\le\|K\|_{L^{1}}\|g_{t}^{c}\|_{C^{r}}
\]
for any $r>\Delta\alpha$. Consequently, $(n\Delta)^{-2}E_{s}\lesssim(n\Delta)^{-2}\E[Y_{1}^{4}]\lesssim n^{-2}\Delta^{-1}$
because $Y_{1}$ has finite fourth moments due to $\|x^{4}\nu\|_{L^{1}}<\infty$.

It remains to bound $E_{s}$ from (\ref{eq:VarM}). Using again $\gamma_0=0$ by Lemma~\ref{lem:drift}, we infer from
$\F[ixP_{\Delta}]=\phi_{\Delta}'=\Delta\psi'\phi=\Delta\F[ix\nu]\phi_{\Delta}$
that
\begin{equation}
xP_{\Delta}=\Delta(x\nu)\ast P_{\Delta}\label{eq:xP}
\end{equation}
and thus $xP_{\Delta}$ has a bounded density satisfying $\|xP_{\Delta}\|_{\infty}\le\Delta\|x\nu\|_{\infty}$.
Together with the Cauchy--Schwarz inequality and Plancherel's identity
we obtain
\begin{align}
E_{s}\le & \Delta\|x\nu\|_{\infty}\int|y|^{3}\Big(\F^{-1}\big[\F g_{t}^{s}(-\bull)m_{\Delta,h}\big](-y)\Big)^{2}\d y\nonumber \\
= & \Delta\|x\nu\|_{\infty}\int\Big|\F^{-1}\big[(\F g_{t}^{s}(-\bull)m_{\Delta,h})''\big](-y)\F^{-1}\big[(\F g_{t}^{s}(-\bull)m_{\Delta,h})'\big](-y)\Big|\d y\nonumber \\
\le & \Delta\|x\nu\|_{\infty}\left\Vert (\F g_{t}^{s}(-\bull)m_{\Delta,h})''\right\Vert _{L^{2}}\left\Vert (\F g_{t}^{s}(-\bull)m_{\Delta,h})'\right\Vert _{L^{2}}.\label{eq:varEs}
\end{align}
The derivatives of the regularized Fourier multiplier are given by
\begin{align}
m_{\Delta,h}'(u) & =ih\frac{\F[xK](hu)}{\phi_{\Delta}(u)}-\Delta\psi'(u)m_{\Delta,h}(u)=ih\frac{\F[xK](hu)}{\phi_{\Delta}(u)}-i\Delta\F[x\nu](u)m_{\Delta,h}(u),\label{eq:mP}\\
m_{\Delta,h}''(u) & =-h^{2}\frac{\F[x^{2}K](hu)}{\phi_{\Delta}(u)}-2i\Delta h\psi'(u)\frac{\F[xK](hu)}{\phi_{\Delta}(u)}-\Delta\big(\psi''(u)(u)+\Delta\psi'(u)^{2}\big)m_{\Delta,h}(u).\label{eq:mPP}
\end{align}
To bound the $L^{2}$-norms in \eqref{eq:varEs}, we use the properties
of $K$, the decay assumption on $\phi_{\Delta}$ and $\psi',\psi''\in L^\infty(\R)$ to obtain 
\begin{align}
\|(\F g_{t}^{s}(-\bull)m_{\Delta,h})'\|_{L^{2}}\le & \|\F[xg_{t}^{s}](-\bull)m_{\Delta,h}\|_{L^{2}}+\|\F[g_{t}^{s}](-\bull)m_{\Delta,h}'\|_{L^{2}}\nonumber \\
\lesssim & (t^{-2}\vee1)(1+h+\Delta)\|(1+|u|)^{\Delta\alpha-1}\|_{L^{2}(I_h)},\nonumber \\
\|(\F g_{t}^{s}(-\bull)m_{\Delta,h})''\|_{L^{2}}\le & \|\F[x^{2}g_{t}^{s}](-\bull)m_{\Delta,h}\|_{L^{2}} +2\|\F[xg_{t}^{s}](-\bull)m_{\Delta,h}'\|_{L^{2}}+\|\F g_{t}^{s}(-\bull)m_{\Delta,h}''\|_{L^{2}}\notag\\
\lesssim & (t^{-2}\vee1)(1+h+\Delta+h^{2}+\Delta h)\|(1+|u|)^{\Delta\alpha-1}\|_{L^{2}(I_h)}.\label{eq:L2gm}
\end{align}
Therefore, $(\Delta n)^{-2}E_{s}\lesssim (t^{-2}\vee1)n^{-2}\Delta^{-1}\|(1+|u|)^{\Delta\alpha-1}\|_{L^{2}(I_h)}^{2}$ which implies
\begin{equation}
\Var(\xi_{k}(t))\lesssim\frac{t^{-2}\vee1}{n^{2}\Delta}\|(1+|u|)^{\Delta\alpha-1}\|_{L^{2}(I_h)}^{2}.\label{eq:VarXi}
\end{equation}
Using (\ref{eq:mP}), (\ref{eq:mPP}), $x\nu\in L^{2}(\R)$, $|\F g_{t}(u)|\lesssim(t^{-2}\vee1)(1+|u|)^{-1}$
and $\|xg_{t}\|_{L^{2}}\lesssim|t|^{-1}$, we deterministically bound
$\xi_{k}(t)$ by 
\begin{align*}
|\xi_{k}(t)|\le & (n\Delta)^{-1}\big\|\F^{-1}[(\F g_{t}(-\bull)m_{\Delta,h})'']\big\|_{\infty}\allowdisplaybreaks\\
\le & (n\Delta)^{-1}\Big(\big\|\F^{-1}\big[\F[x^{2}g_{t}](-\bull)m_{\Delta,h}\big]\big\|_{\infty}+2\|\F[xg_{t}](-\bull)m_{\Delta,h}'\|_{L^{1}}+\|\F g_{t}(-\bull)m_{\Delta,h}''\|_{L^{1}}\Big)\allowdisplaybreaks\\
\lesssim & (n\Delta)^{-1}\Big(\big\|\F^{-1}\big[\F[x^{2}g_{t}](-\bull)m_{\Delta,h}\big]\big\|_{\infty}+\|xg_{t}\|_{L^{2}}\|m_{\Delta,h}'\|_{L^{2}}\\
 & +\|(1+|u|)\F g_{t}(u)\|_{\infty}\|(1+|u|)^{-1}m_{\Delta,h}''(u)\|_{L^{1}}\Big)\allowdisplaybreaks\\
\lesssim & (n\Delta)^{-1}\Big(\|\F[x^{2}g_{t}^{s}](-\bull)m_{\Delta,h}\|_{L^{1}}+\|\F^{-1}\big[\F[x^{2}g_{t}^{c}](-\bull)m_{\Delta,h}\big]\|_{\infty}\\
 & +\|xg_{t}\|_{L^{2}}\big(h\|(1+|u|)^{\Delta\alpha}\|_{L^{2}(I_h)}+\Delta \|x\nu\|_{L^{1}}\|(1+|u|)^{\Delta\alpha}\|_{L^{2}(I_h)}\big)\\
 & +(t^{-2}\vee1)(\Delta+h^{2}+\Delta h)\|(1+|u|)^{\Delta\alpha-1}\|_{L^{1}(I_h)}\Big)\allowdisplaybreaks\\
\lesssim & (n\Delta)^{-1}\|\F[x^{2}g_{t}^{s}](-\bull)m_{\Delta,h}\|_{L^{1}}+(n\Delta)^{-1}\|\F^{-1}\big[\F[x^{2}g_{t}^{c}](-\bull)m_{\Delta,h}\big]\|_{\infty}\\
 & +(n\Delta)^{-1}(t^{-2}\vee1)h^{-1/2}\|(1+|u|)^{\Delta\alpha-1}\|_{L^{2}(I_h)}.
\end{align*}
The term corresponding to the singular part $g_{t}^{s}$ in the previous bound can be estimated by
\begin{align*}
(n\Delta)^{-1}\big\|\F[x^{2}g_{t}^{s}](-\bull)m_{\Delta,h}\big\|_{L^{1}} & \lesssim(n\Delta)^{-1}\|(1+|u|)\F[x^{2}g_{t}^{s}](u)\|_{\infty}\|(1+|u|)^{\Delta\alpha-1}\|_{L^{1}(I_h)}\\
 & \lesssim(n\Delta)^{-1}(t^{-2}\vee1)h^{-1/2}\|(1+|u|)^{\Delta\alpha-1}\|_{L^{2}(I_h)}.
\end{align*}
For the continuous part $g_{t}^{c}$ we apply the Fourier multiplier
theorem as above to see that $$\|\F^{-1}\big[\F[x^{2}g_{t}^{c}](-\bull)m_{\Delta,h}\big]\|_{\infty}\lesssim\|x^{2}g_{t}^{c}\|_{C^{r}}$$
for any $r>\Delta\alpha$. Therefore,
\begin{align}
|\xi_{k}(t)|\le & (n\Delta)^{-1}\big\|\F^{-1}[(\F g_{t}(-\bull)m_{\Delta,h})'']\big\|_{\infty}\nonumber \\
\lesssim & (n\Delta)^{-1}(t^{-2}\vee1)h^{-1/2}\|(1+|u|)^{\Delta\alpha-1}\|_{L^{2}(I_h)}.\label{eq:DetXi}
\end{align}
Using (\ref{eq:VarXi}) and (\ref{eq:DetXi}), Bernstein's inequality
yields for some constant $c>0$ the claimed concentration result.
\end{proof}

Combining the previous results, we obtain minimax convergence rates
for estimating the (generalized) distribution function of the jump
measure.
\begin{proof}[Proof of Proposition \ref{prop:DistFunct}]
In the following, $t$ is fixed and thus omitted in the constants.
Using the error decomposition (\ref{eq:Decomp}), Lemma~\ref{lem:vola},
Proposition~\ref{prop:bias}, we obtain
\begin{align*}
|\hat{N}_{h}(t)-N(t)|\le & |B_{n}(t)|+|S_{n}(t)|+|V_{n}(t)|\\
\lesssim & h^{s+1}+|S_{n}(t)|.
\end{align*}
Using $|\F g_{t}(u)|\lesssim(1+|u|)^{-1}$, Plancherel's identity
and Lemma~\ref{lem:Remainder} yield for the stochastic error from
(\ref{eq:Decomp}) and the linearized stochastic error term $L_{\Delta,n}$
defined in (\ref{eq:linError}) 
\begin{align}
 & \big|S_{n}(t)-L_{\Delta,n}(t)\big|\nonumber \\
= & \left|\int g_{t}(x)\F^{-1}\Big[\F K(hu)\Big(\hat{\psi}''_{n}(u)-\psi''(u)-\Delta^{-1}\big(\phi_{\Delta}^{-1}(\phi_{\Delta,n}-\phi_{\Delta})\big)''(u)\Big)\Big](x)\d x\right|\nonumber \\
\lesssim & \sup_{\left|u\right|\le h^{-1}}\big|\hat{\psi}_{n}''(u)-\psi''(u)-\Delta^{-1}\big(\phi_{\Delta}^{-1}(\phi_{\Delta,n}-\phi_{\Delta})\big)''(u)\big|\int(1+|u|)^{-1}|\F K(hu)|\d u\nonumber \\
\le & \sup_{\left|u\right|\le h^{-1}}\big|\hat{\psi}_{n}''(u)-\psi''(u)-\Delta^{-1}\big(\phi_{\Delta}^{-1}(\phi_{\Delta,n}-\phi_{\Delta})\big)''(u)\big|\|K\|_{L^{1}}\int_{-1/h}^{1/h}(1+|u|)^{-1}\d u\nonumber \\
= & \mathcal{O}_{P}\big((\Delta\|\psi'\big\|_{L^{\infty}(I_h)}+\Delta^{3/2}\|\psi'\big\|_{L^{\infty}(I_h)}^{2}+1)\frac{|\log h|^{2+\delta}}{n\Delta^{1/2}}\|\phi_{\Delta}^{-1}\|_{L^{\infty}(I_h)}^{2}\big),\label{eq:linearization}
\end{align}
provided that $n^{-1/2}(\log h^{-1})^{(1+\delta)/2}\|\phi_{\Delta}^{-1}\|_{L^{\infty}(I_h)}\to0$.
The latter condition is satisfied for the choices $h=h_{n,\Delta}$ in
both cases. 

Let $(\sigma^{2},\gamma,\nu)\in\mathcal{D}^{s}(\alpha,m,U,R)$. We conclude
from (\ref{eq:linearization}), where $\psi'$ is uniformly bounded
and $\phi_{\Delta}$ decays polynomially, and Propositions~\ref{prop:variance}
and \ref{prop:mainStochErr} for $\Delta\alpha>1/2$ 
\begin{align*}
S_{n}(t)= & L_{\Delta,n}(t)+\mathcal{O}_{P}\big(n^{-1}\Delta^{-1/2}\log(h^{-1})^{2+\delta}h^{-2\Delta\alpha}\big)\\
= & \mathcal{O}_{P}\big((n\Delta)^{-1/2}h^{-\Delta\alpha+1/2}+n^{-1}\Delta^{-1/2}\log(h^{-1})^{2+\delta}h^{-2\Delta\alpha}\big)\\
= & \mathcal{O}_{P}\big((n\Delta)^{-1/2}h^{-\Delta\alpha+1/2}(1+n^{-1/2}\log(h^{-1})^{2+\delta}h^{-\Delta\alpha-1/2})\big).
\end{align*}
Therefore, we obtain the rate $r_{n,\Delta}$ by plugging in $h=h_{n,\Delta}=(n\Delta)^{-1/(2s+2\Delta\alpha+1)}$ and similarly for $\Delta\alpha\le 1/2$.

Let us consider the case $(\sigma^{2},\gamma,\nu)\in\mathcal{E}^{s}(\beta,m,U,r,R)$.
Owing to $\|\psi'\big\|_{L^{\infty}(I_h)}\lesssim h^{-1}$ and
the exponential decay of $\phi_{\Delta}$, we obtain from (\ref{eq:linearization})
and Proposition~\ref{prop:variance} that 
\begin{align*}
S_{n}(t)= & \mathcal{O}_{P}\big((n\Delta)^{-1/2}\big(\log(h^{-1})+\Delta h^{-1}+\Delta^{3/2}h^{-2}\\
 & \qquad+n^{-1/2}\log(h^{-1})^{2+\delta}(\Delta h^{-1}+\Delta^{3/2}h^{-2}+1)\big)\exp(r\Delta h^{-\beta})\big).
\end{align*}
Therefore, plugging in $h=h_{n,\Delta}=(\frac{\delta}{2r}\frac{\log(n\Delta)}{\Delta})^{-1/\beta},\delta\in(0,1)$,
yields
\begin{align*}
&\big|\hat{N}_{h}(t)-N(t)\big|\\
&\quad=\mathcal{O}_{P}\Big(\Big(\frac{\log(n\Delta)}{\Delta}\Big)^{-(s+1)/\beta}+(n\Delta)^{-(1-\delta)/2}\log(n\Delta)^{2/\beta}(|\log\Delta|+\Delta^{1-1/\beta}+\Delta^{3/2-2/\beta})\Big).
\end{align*}

\end{proof}

\subsubsection{Uniform loss for density estimation}
Applying a similar decomposition as in (\ref{eq:Decomp}) and the linearization
Lemma~\ref{lem:Remainder}, we obtain
\begin{align}
  &\hat{\nu}_{h}(t)-\nu(t)\notag\\
  =&\frac{1}{t^2}\left(\big(K_{h}\ast\big(y^{2}\nu\big)-x^{2}\nu\big)(t)+\F^{-1}\Big[\Big(\psi''(u)-\hat{\psi}''_{n}(u)\Big)\F K(hu)\Big](t)+\sigma^{2}K_{h}(t)\right)\label{eq:DecompDens}\\
= & \frac{1}{t^2}\Big(\underbrace{(K_{h}\ast(y^{2}\nu)-x^{2}\nu)(t)}_{=:B^{\nu}(t)}-\underbrace{\frac{1}{\Delta}\F^{-1}\Big[\F K(hu)\Big(\frac{\phi_{\Delta,n}-\phi_{\Delta}}{\phi_{\Delta}}\Big)''(u)\Big](t)}_{=:L^\nu_{\Delta,n}(t)}+R_{\Delta,n}+\sigma^{2}K_{h}(t)\Big)\notag
\end{align}
for some remainder $R_{\Delta,n}$ which is of order
\begin{align}
|R_{\Delta,n}| & =\mathcal{O}_{P}\big((\Delta\|\psi'\big\|_{L^{\infty}(I_h)}+\Delta^{3/2}\|\psi'\big\|_{L^{\infty}(I_h)}^{2}+1)\label{eq:RemainderDens}\\
 & \qquad\qquad\times n^{-1}\Delta^{-1/2}h^{-1}\log(h^{-1})^{1+\delta}\|\phi_{\Delta}^{-1}\|_{L^{\infty}(I_h)}^{2}\big).\nonumber 
\end{align}
We will need the following concentration result for the main stochastic
error term of the density estimation problem
\begin{equation}
M^\nu_{\Delta,n}(t)=-\frac{1}{n\Delta}\sum_{k=1}^{n}\F^{-1}\Big[m_{\Delta,h}\big(Y_{k}^{2}e^{iuY_{k}}-\E[Y_{k}^{2}e^{iuY_{k}}]\big)\Big](t),\label{eq:MainStochErr}
\end{equation}
where we recall $m_{\Delta,h}=\F K(h\bull)/\phi_{\Delta}$. We will prove it analogously to Proposition~\ref{prop:mainStochErr}.
\begin{lemma}\label{lem:concMain} 
  Let $\alpha,R>0,m>4,U\subset\R$ and $t\neq0$ and let the kernel satisfy \eqref{eq:propKernel} for $p\ge1$. If $(\sigma^{2},\gamma,\nu)\in\mathcal{D}^{s}(\alpha,m,U,R)$, then there is some constant $c>0$, depending only on $\alpha,R$,
  such that for any $\kappa_{0}>0$ and any $n\in\mathbb{N},\Delta,h>0$
  \[
  P\left(|M^\nu_{\Delta,n}(t)|>\kappa_{0}(\Delta n)^{-1/2}h^{-\Delta\alpha-1/2}\right)\le2\exp\left(-\frac{c\kappa_{0}^{2}}{(1+t^{3})(1+\kappa_{0}(n\Delta h)^{-1/2})}\right).
  \]
\end{lemma}
\begin{proof}
We apply Bernstein's inequality to the sum of the independent and
centered random variables
\[
M^\nu_{\Delta,n}(t)=-\sum_{k=1}^{n}\big(\xi_{k}-\E[\xi_{k}]\big),\quad\text{with}\quad\xi_{k}:=\frac{1}{n\Delta}\F^{-1}\Big[m_{\Delta,h}(u)Y_{k}^{2}e^{iuY_{k}}\Big](t).
\]
$\Var(\xi_{k})$ can be estimated similarly to (\ref{eq:varEs}).
We obtain by (\ref{eq:xP}), the Cauchy--Schwarz inequality, Plancherel's
identity, \eqref{eq:mP} and \eqref{eq:mPP} that
\begin{align*}
\Var(\xi_{k})\le&\E\big[\xi_{k}^{2}\big]=  (n\Delta)^{-2}\E\Big[Y_{k}^{4}\Big(\F^{-1}\Big[m_{\Delta,h}(u)e^{-iut}\Big](-Y_{k})\Big)^{2}\Big]\\
\le & \frac{1}{n^{2}\Delta}\|x\nu\|_{\infty}\Big\|y^{2}\F^{-1}\Big[m_{\Delta,h}(u)e^{-iut}\Big](-y)\Big\|_{L^{2}}\Big\|y\F^{-1}\Big[m_{\Delta,h}(u)e^{-iut}\Big](-y)\Big\|_{L^{2}}\\
= & \frac{1}{n^{2}\Delta}\frac{\|x\nu\|_{\infty}}{2\pi}\|(m_{\Delta,h}(u)e^{-iut})''\|_{L^{\text{2}}}\|(m_{\Delta,h}(u)e^{-iut})'\|_{L^{2}}\\
\lesssim & \frac{1}{n^{2}\Delta}(1+t^{3})\|(1+|u|)^{\Delta\alpha}\|_{L^{2}(I_h)}^{2}\lesssim\frac{1}{n^{2}\Delta}(1+t^{3})h^{-2\Delta\alpha-1}.
\end{align*}
Moreover, $\xi_{k}$ admits the deterministic bound
\begin{align*}
|\xi_{k}|= & \frac{1}{n\Delta}Y_{k}^{2}\big|\F^{-1}\Big[m_{\Delta,h}(u)e^{-iut}\Big](-Y_{k})\big|\\
\le & \frac{1}{n\Delta}\big|\F^{-1}\big[(m_{\Delta,h}(u)e^{-iut})''\big](-Y_{k})\big|\\
\le & \frac{1}{2\pi n\Delta}\left(\|m_{\Delta,h}''\|_{L^{1}}+2t\|m_{\Delta,h}'\|_{L^1}+t^{\text{2}}\|m_{\Delta,h}\|_{L^{1}}\right)\\
\lesssim & \frac{1}{n\Delta}(1+t^{2})\|(1+|u|)^{\Delta\alpha}\|_{L^{1}(I_h)}\lesssim\frac{1+t^{2}}{n\Delta}h^{-\Delta\alpha-1}.
\end{align*}
Therefore, Bernstein's inequality yields for a constant $c>0$ and
any $\kappa>0$ 
\[
P\left(|M^\nu_{\Delta,n}(t)|>\kappa\right)\le2\exp\left(-\frac{c\Delta n\kappa^{2}}{(1+t^{3})h^{-2\Delta\alpha-1}+\kappa(1+t^{2})h^{-\Delta\alpha-1}}\right).
\]
Choosing $\kappa=\kappa_{0}(\Delta n)^{-1/2}h^{-\Delta\alpha-1/2}$
for $\kappa_{0}>0$, we conclude
\begin{align*}
P\left(|M^\nu_{\Delta,n}(t)|>\kappa_{0}(\Delta n)^{-1/2}h^{-\Delta\alpha-1/2}\right)\le & 2\exp\left(-\frac{c\kappa_{0}^{2}}{(1+t^{3})(1+\kappa_{0}(n\Delta h)^{-1/2})}\right).\qedhere
\end{align*}

\end{proof}

Proposition~\ref{prop:DensLevyRate} is an immediate consequence of the following result.
\begin{prop}\label{prop:UniformDens} 
  Let $\alpha,\beta,s,r,R>0,m>4$, let the kernel satisfy (\ref{eq:propKernel})
  with order $p\ge s$ and let $U\subset\R$ be a bounded, open set
  which is bounded away from zero. Then we have
  \begin{enumerate}
  \item uniformly in $(\sigma^{2},\gamma,\nu)\in\mathcal{D}^{s}(\alpha,m,U,R)$,
  if $\log(n\Delta)/(n\Delta h)\to0$ , 
  \[
  \sup_{t\in U}|\hat{\nu}_{h}(t)-\nu(t)|=\mathcal{O}_{P,\mathcal{D}^{s}}\Big(h^{s}+\Big(\frac{\log n\Delta}{n\Delta}\Big)^{1/2}h^{-\Delta\alpha-1/2}\Big),
  \]

  \item uniformly in \textup{$(\sigma^{2},\gamma,\nu)\in\mathcal{E}^{s}(\beta,m,U,r,R)$}
  \[
  \sup_{t\in U}|\hat{\nu}_{h}(t)-\nu(t)|=\mathcal{O}_{P,\mathcal{E}^{s}}\Big(h^{s}+(n\Delta)^{-1/2}(h^{-1}+\Delta h^{-2}+\Delta^{3/2}h^{-3})e^{r\Delta h^{-\beta}}\Big).
  \]

  \end{enumerate}
\end{prop}
\begin{proof}
We start with the error decomposition \eqref{eq:DecompDens}. By standard approximation arguments the deterministic error satisfies
$\sup_{t\in U}|B^{\nu}(t)|\lesssim h^{s}$ if $\nu\in C^{s}(U)$ for
an open set $U$ , $s>0$ and if the kernel satisfies (\ref{eq:propKernel})
with order $p\ge s$. Moreover, the assumptions $\supp\F K\subset[-1,1]$
and $x^{p+1}K(x)\in L^{1}(\R)$ imply $\|x^{p+1}K(x)\|_{\infty}<\infty$
which yields 
\[
|\sigma^{2}K_{h}(t)|\le\sigma^{2}h^{-1}\sup_{|x|>|t|/h}|K(x)|\lesssim\sigma^{2}|t|^{-s-1}h^{s}.
\]
Since $U$ is bounded away from zero, the previous display gives a
uniform bound on $U$. Using the main stochastic error term $M^\nu_{\Delta,n}$
from (\ref{eq:MainStochErr}), the linearized stochastic error term from \eqref{eq:DecompDens}
can be decomposed similarly to (\ref{eq:DecompLin}) into
\begin{align}
L^\nu_{\Delta,n}(t)=&M^\nu_{\Delta,n}(t)+2\F^{-1}\big[m_{\Delta,h}\psi'(\phi_{\Delta,n}-\phi_{\Delta})'\big](t)\notag\\
&\qquad+{\cal F}^{-1}[m_{\Delta,h}\big(\psi''-\Delta(\psi')^{2}\big)(\phi_{\Delta,n}-\phi_{\Delta})](t).\label{eq:decompLNu}
\end{align}
To derive the appropriate bounds for $R_{\Delta,n}$ and $L^\nu_{\Delta,n}$,
we will distinguish again between the mildly and the severely ill-posed
case.

We start with the severely ill-posed case $(\sigma^{2},\gamma,\nu)\in\mathcal{E}^{s}(\beta,m,U,r,R).$
Using Fubini's theorem and the properties of $m_{\Delta,h}$ as well as $|\psi''(u)|\lesssim1,|\psi'(u)|\lesssim1+|u|,u\in\R,$ and (\ref{eq:phiMom}), we obtain
\begin{align*}
&\E[\sup_{t\in U}|L^\nu_{\Delta,n}(t)|]\\
&\qquad\le  \Delta^{-1}\E\Big[\|{\cal F}^{-1}\Big[m_{\Delta,h}(\phi_{\Delta,n}''-\phi_{\Delta}'')\|_{\infty}+2\E\Big[\|\F^{-1}\big[m_{\Delta,h}\psi'(\phi_{\Delta,n}-\phi_{\Delta})'\big]\|_{\infty}\Big]\\
 &\qquad\qquad +\E\Big[\|{\cal F}^{-1}[m_{\Delta,h}\big(\psi''-\Delta(\psi'){}^{2}\big)(\phi_{\Delta,n}-\phi_{\Delta})]\|_{\infty}\Big]\\
&\qquad\lesssim  \int_{-1/h}^{1/h}\Big(\Delta^{-1}\E[(\phi_{\Delta,n}''(u)-\phi_{\Delta}''(u))^{2}]^{1/2}+\E[(\phi_{\Delta,n}'(u)-\phi_{\Delta}'(u))^{2}]^{1/2}(1+|u|)\\
 & \qquad\qquad+\E[(\phi_{\Delta,n}(u)-\phi_{\Delta}(u))^{2}]^{1/2}(1+\Delta(1+|u|)^{2})\Big)\exp(r\Delta|u|^{\beta})\d u\\
&\qquad\le  (n\Delta)^{-1/2}\int_{-1/h}^{1/h}\Big(1+\Delta(1+|u|)+\Delta^{1/2}+\Delta^{3/2}(1+|u|)^{2}\Big)\exp(r\Delta|u|^{\beta})\d u\\
&\qquad\lesssim  (n\Delta)^{-1/2}(h^{-1}+\Delta h^{-2}+\Delta^{3/2}h^{-3})\exp(r\Delta h^{-\beta}).
\end{align*}
The remainder (\ref{eq:RemainderDens}) is of smaller order 
\[
|R_{\Delta,n}|=\mathcal O_{P}\Big(n^{-1}\Delta^{-1/2}\log(h^{-1})^{1+\delta}(h^{-1}+\Delta h^{-2}+\Delta^{3/2}h^{-3})\big)\exp(2r\Delta h^{-\beta})\Big).
\]

Now let us consider $(\sigma^{2},\gamma,\nu)\in\mathcal{D}^{s}(\alpha,m,U,R)$,
where we have
\[
|R_{\Delta,n}|=\mathcal{O}_{P}\big((\Delta+\Delta^{3/2}+1)n^{-1}\Delta^{-1/2}\log(h^{-1})^{1+\delta}h^{-2\Delta\alpha-1}\big).
\]
To bound $\sup_{t\in U}|L^\nu_{\Delta,n}(t)|$, we note for the second
and the third term in (\ref{eq:decompLNu}) that similarly to (\ref{eq:nonCritMild})
with the Cauchy--Schwarz inequality and Fubini's theorem 
\begin{align*}
 & \E\Big[\sup_{t\in U}|L^\nu_{\Delta,n}(t)-M^\nu_{\Delta,n}(t)|\Big]\\
\le & 2\E\Big[\|m_{\Delta,h}\psi'(\phi_{\Delta,n}-\phi_{\Delta})'\|_{L^{1}}\Big]+\E\Big[\|m_{\Delta,h}\big(\psi''-\Delta(\psi'){}^{2}\big)(\phi_{\Delta,n}-\phi_{\Delta})\|_{L^{1}}\Big]\\
\lesssim & 2\|x\nu\|_{L^{2}}\E\big[\|m_{\Delta,h}(\phi_{\Delta,n}-\phi_{\Delta})'\|_{L^{2}}\big]+\|x^{2}\nu\|_{L^{2}}\E\big[\|m_{\Delta,h}(\phi_{\Delta,n}-\phi_{\Delta})\|_{L^{2}}\big]\\
 & +\Delta\|x\nu\|_{L^{2}}\E\big[\|m_{\Delta,h}\psi'(\phi_{\Delta,n}-\phi_{\Delta})\|_{L^{2}}\big]\\
\lesssim & n^{-1/2}\Big(2\Delta^{1/2}\|x\nu\|_{L^{2}}+\|x^{2}\nu\|_{L^{2}}+\Delta\|x\nu\|_{L^{2}}\|x\nu\|_{L^{1}}\Big)\|m_{\Delta,h}\|_{L^{2}}\\
\lesssim & n^{-1/2}h^{-\Delta\alpha-1/2}.
\end{align*}
It remains to estimate the main stochastic term $M^\nu_{\Delta,n}$.
Since $U$ is bounded, we find a finite number of points $t_{1},\dots,t_{L_{n}}\in U$
such that $\sup_{t\in U}\min_{l=1,\dots,L_{n}}|t-t_{l}|\le(n\Delta)^{-2}$
and $L_{n}\lesssim(n\Delta)^{2}$. The inequalities by Young and by Cauchy--Schwarz together with $(n\Delta)^{-1}\lesssim1$ yield
\begin{align*}
 & \E\Big[\sup_{t\in U}\min_{l=1,\dots,L_{n}}|M^\nu_{\Delta,n}(t)-M^\nu_{\Delta,n}(t_{l})|\Big]\\
\lesssim & (n\Delta)^{-2}\Delta^{-1}\E\Big[\big\|(K_{h}')\ast\big(\F\big[\1_{[-1/h,1/h]}\phi_{\Delta}^{-1}(\phi_{\Delta,n}''-\phi_{\Delta}'')\big]\big)\big\|_{\infty}\Big]\\
\lesssim & (n\Delta)^{-2}\Delta^{-1}h^{-1}\|K'\|_{L^{1}}\int_{-1/h}^{1/h}|\phi_{\Delta}^{-1}|(u)\E[|\phi_{\Delta,n}''-\phi_{\Delta}''|(u)]\d u\\
\lesssim & (n\Delta)^{-2}h^{-\Delta\alpha-2}.
\end{align*}
Together with Markov's inequality and Lemma~\ref{lem:concMain}, defining $T:=\sup_{t\in U}{|t|}$ this
yields for $n\Delta$ sufficiently large 
\begin{align*}
 & P\Big(\sup_{t\in U}|M^\nu_{\Delta,n}(t)|>\kappa_{0}\Big(\frac{\log(n\Delta)}{n\Delta}\Big)^{1/2}h^{-\Delta\alpha-1/2}\Big)\\
\le & P\Big(\max_{l=1,\dots,L_{n}}|M^\nu_{\Delta,n}(t_{k})|>\frac{\kappa_{0}}{2}\Big(\frac{\log(n\Delta)}{n\Delta}\Big)^{1/2}h^{-\Delta\alpha-1/2}\Big)\\
 & +\frac{2}{\kappa_{0}}\big(\frac{n\Delta}{\log(n\Delta)}\big)^{1/2}h^{\Delta\alpha+1/2}\E\Big[\sup_{t\in U}\min_{l=1,\dots,L_{n}}|M^\nu_{\Delta,n}(t)-M^\nu_{\Delta,n}(t_{l})|\Big]\\
\le & 2L_{n}\exp\Big(-\frac{c\kappa_{0}^{2}\log(n\Delta)}{2(1+T^3)(2+\kappa_{0}(\log(n\Delta)/(n\Delta h))^{1/2})}\Big)+o\Big((n\Delta)^{-3/2}(\log n\Delta)^{-1/2}h^{-3/2}\Big)\\
\le & 2\exp\big((2-\tfrac{c}{6}\kappa_{0}^{2}(1+T^3)^{-1})\log(n\Delta)\big)+o(1),
\end{align*}
which converges to zero as $n\Delta\to\infty$ if $\kappa_{0}$ is
chosen sufficiently large.
\end{proof}

\subsubsection{Proof of Theorem~\ref{thm:rateQuantiles}}

In view of the error representation \eqref{eq:errorRep}, the missing ingredient to prove Theorem~\ref{thm:rateQuantiles} is consistency of $\hat q_{\tau,h}^\pm$. We prove it similarly to \cite{dattnerEtAl2013}, but since $\nu$ has no bounded density and we minimize on an unbounded interval the lemma is more involved. 
\begin{lemma}
  \label{lem:consistency}Let $\alpha,\beta,s,\zeta,r,R>0,m>4,s'\in(-1,0]$
  and let $\eta_{n}\downarrow0$ with $\eta_{n}^{-1}\lesssim\log n$.
  Suppose the kernel satisfies (\ref{eq:propKernel}) with order $p\ge1$. Then we have
  \begin{enumerate}
  \item for any bandwidth satisfying $(\log n)^{4}h^{1+s'}\to0$ and $(n\Delta)^{-1}h^{-2\Delta\alpha-1}\to0$
  \[
  \sup_{(\sigma^{2},\gamma,\nu)\in\mathcal{\tilde{D}}_{\tau}^{s,s'}(\alpha,m,\zeta,\eta_{n},R)}P\big(|\hat{q}_{\tau,h}^{\pm}-q_{\tau}^{\pm}|>\delta\big)\to0\quad\text{for all }\delta>0,
  \]

  \item for any bandwidth satisfying $(\log n)^{4}h^{1+s'}\to0$ and $(\log n)^{4}(n\Delta)^{-1}h^{-2}e^{2r\Delta h^{-\beta}}\to0$
  \[
  \sup_{(\sigma^{2},\gamma,\nu)\in\mathcal{\tilde{E}}_{\tau}^{s,s'}(\beta,m,\zeta,\eta_{n},r,R)}P\big(|\hat{q}_{\tau,h}^{\pm}-q_{\tau}^{\pm}|>\delta\big)\to0\quad\text{for all }\delta>0.
  \]
  \end{enumerate}
\end{lemma}
\begin{proof}
We adopt the general strategy of the proof of Theorem 5.7 by \citet{vanderVaart1998}
in the classical M-estimation setting. Without loss of generality,
we only consider $\hat{q}_{\tau,h}^{+}$.

\emph{Step 1:} By the H\"older regularity we have $\nu(t)\ge\nu(q_{\tau}^{+})-|\nu(q_{\tau}^{+})-\nu(t)|\ge\frac{1}{R}-R|q_{\tau}^{+}-t|^{1\wedge\alpha}\ge\frac{1}{2R}$
for $|q_{\tau}^{+}-t|\le(2R^{2})^{-(1\vee\alpha^{-1})}$. Without
loss of generality we can assume $\delta\le\delta_0:=(2R^{2})^{-(1\vee\alpha^{-1})}\wedge(\eta_n/2)\wedge\zeta$,
otherwise consider $\delta\wedge\delta_0$.
Using $N(q_{\tau}^{+})=\tau$ and monotonicity of $N$, we obtain
the uniqueness condition 
\begin{align}\label{eq:Unique}
\inf_{t>\eta_{n}:|t-q_{\tau}^{+}|\ge\delta}|N(t)-\tau| & \ge\delta\inf_{t>\eta_{n}:|t-q_{\tau}^{+}|\le\delta}\nu(t)\ge\frac{\delta}{2R}.
\end{align}

\emph{Step 2:} We construct an event $A$ with $P(A)\to1$ such that
  \begin{equation}\label{eq:ZEst}
    \hat N_h(\hat q^+_\tau)-\tau=0\quad\text{on }A.
  \end{equation}
  Using $N(q^+_\tau)=\tau$, monotonicity of $N$ and \eqref{eq:Unique}, we conclude for $\delta\in(0,\delta_0)$
  \[
    N(q^+_\tau+\delta)-\tau\le-\frac{\delta }{2R}<0<\frac{\delta }{2R}\le N(q^+_\tau-\delta)-\tau.
  \]
  Proposition~\ref{prop:DistFunct} implies that
  \begin{equation}\label{eq:A5}
    A:=\big\{\forall \eta\in\{q_\tau-\delta,q_\tau+\delta\}:|\hat N_h(\eta)-N(\eta)|\le\frac{\delta }{4R}\big\}\quad\text{satisfies}\quad P(A)\to1.
  \end{equation}
  We conclude on $A$ that $\hat N_h(q_\tau+\delta)<\tau< \hat N_h(q_\tau-\delta)$. Since $\hat N_h$ admits a derivative, namely $\hat\nu_h$, it is continuous and thus there is an intermediate point $\xi_h\in(q^+_\tau-\delta,q^+_\tau+\delta)$ such that $\hat N_h(\xi_h)=\tau$. Since $\hat{q}_{\tau,h}^{+}$ minimizes $|\hat{N}_{h}(\bull)-\tau|$
on the interval $(\eta_{n},\infty)$ and $q_{\tau}-\delta\in(\eta_{n},\infty)$ for $\eta_n$ and $\delta$ sufficiently small, we obtain \eqref{eq:ZEst}.

\emph{Step 3:} We infer from Steps~1 and 2 
\begin{align}
P\big(|\hat{q}_{\tau,h}^{+}-q_{\tau}^{+}|>\delta\big) & \le P\big(|N(\hat{q}_{\tau,h}^{+})-\tau|\ge\delta/(2R)\big)\nonumber \\
 & =P\big(|N(\hat{q}_{\tau,h}^{+})-\hat{N}_{h}(\hat{q}_{\tau,h}^{+})|\ge\delta/(2R)\big)+o(1)\nonumber \\
 & \le P\big(\sup_{t\in(\eta_{n},\infty)}|N(t)-\hat{N}_{h}(t)|\ge\delta/(2R)\big)+o(1).\label{eq:ConsistencyBoundLevy}
\end{align}
Hence, it remains to show uniform consistency of $\hat{N}_{h}(t)$.
Applying the error decomposition~(\ref{eq:Decomp}) $|N(t)-\hat{N}_{h}(t)|\le|B_{n}(t)|+|S_{n}(t)|+|V_{n}(t)|$
and the estimates in Lemma~\ref{lem:vola} and Proposition~\ref{prop:bias}
we obtain
\begin{align*}
\sup_{t\in(\eta_{n},\infty)}|V_{n}(t)|\lesssim & \eta_{n}^{-4}h,\qquad\sup_{t\in(\eta_{n},\infty)}|B_{n}(t)|\lesssim\eta_{n}^{-4}h^{1+s'}.
\end{align*}
According to Lemma~\ref{lem:Remainder} and Proposition~\ref{prop:variance},
the stochastic error term can be decomposed into $S_{n}(t)=M_{\Delta,n}(t)+(L_{\Delta,n}-M_{\Delta,n})(t)+R_{n}(t),$
where $\sup_{|t|\ge\eta_{n}}|L_{\Delta,n}-M_{\Delta,n}|(t)+|R_{n}|(t)$
is of the order claimed in Proposition~\ref{prop:variance}. For
$(\sigma^{2},\gamma,\nu)\in\mathcal{E}^{s}(\beta,m,U,r,R)$ the main
stochastic error term is uniformly bounded as well. For the case $(\sigma^{2},\gamma,\nu)\in\mathcal{D}^{s}(\alpha,m,U,R)$
it remains to apply Proposition~\ref{prop:mainStochErr} on an appropriate
grid. For $\kappa\in(0,1)$ we define a grid $v_{l}=-\kappa^{-1}+l\kappa$
with $l=0,\dots,L:=2\lfloor\kappa^{-2}\rfloor$. Set $t_l:=\sign (v_l)(|v_l|\vee\eta_n)$ for $l=0,\dots,L$. Then
\begin{align*}
\sup_{|t|\in(\eta_{n},\infty)}|M_{\Delta,n}(t)|\le & \sup_{l=0,\dots,L}|M_{\Delta,n}(t_{l})|+\sup_{t\in(\eta_{n},\infty)}\min_{l=0,\dots,L}|M_{\Delta,n}(t)-M_{\Delta,n}(t_{l})|.
\end{align*}
Noting that $\|g_{t}-g_{s}\|_{L^{1}}\le|t-s|/(s^{2}\wedge t^{2})$
and $\|g_{t}-g_{s}\|_{L^{1}}\le2\int_{|t|\wedge|s|}^{\infty}x^{-2}\d x\sim|t|^{-1}\vee|s|^{-1}$,
increments of $M_{\Delta,n}(t)$ can be estimated using Plancherel's
identity and Fubini's theorem
\begin{align*}
 & \E\Big[\sup_{t\in(\eta_{n},\infty)}\min_{l=0,\dots,L}|M_{\Delta,n}(t)-M_{\Delta,n}(t_{l})|\Big]\\
\le & \frac{1}{2\pi\Delta}\E\Big[\sup_{t\in(\eta_{n},\infty)}\min_{l=0,\dots,L}\Big|\int\F[g_{t}-g_{t_{l}}](-u)m_{\Delta,h}(u)(\phi_{\Delta,n}''(u)-\phi_{\Delta}''(u))\d u\Big|\Big]\\
\le & \frac{1}{2\pi\Delta}\sup_{t\in(\eta_{n},\infty)}\min_{l=0,\dots,L}\|g_{t}-g_{t_{l}}\|_{L^{1}}\int_{-1/h}^{1/h}|m_{\Delta,h}(u)|\E[(\phi_{\Delta,n}''(u)-\phi_{\Delta}''(u))^{2}]^{1/2}\d u\\
\lesssim & (n\Delta)^{-1/2}\sup_{t\in(\eta_{n},\infty)}\min_{l=0,\dots,L}\|g_{t}-g_{t_{l}}\|_{L^{1}}\|(1+|u|)^{\Delta\alpha}\|_{L^{1}(I_h)}\\
\lesssim & (n\Delta)^{-1/2}\eta_{n}^{-2}\kappa h^{-\Delta\alpha-1}.
\end{align*}
Choosing $\kappa=\eta_{n}^{2}(n\Delta)^{-1/2}$, Markov's inequality
and Proposition~\ref{prop:mainStochErr} yield for $\Delta n$ sufficiently
large and some constant $c>0$
\begin{align*}
 & P\Big(\sup_{|t|\in(\eta_{n},\infty)}|M_{\Delta,n}(t)|>\delta\Big)\\
\le & P\Big(\sup_{l=0,\dots,L}|M_{\Delta,n}(t_{j})|>\frac{\delta}{2}\Big)+\frac{2}{\delta}\E\Big[\sup_{t\in(\eta_{n},\infty)}\min_{l=0,\dots,L}|M_{\Delta,n}(t)-M_{\Delta,n}(t_{l})|\Big]\\
\le & 2(L+1)\exp\Big(-c\delta^{2}\eta_{n}^{3}n\Delta\|(1+|u|)^{\Delta\alpha-1}\|_{L^{2}(I_h)}^{-2}\Big)+(n\Delta)^{-1/2}\eta_{n}^{-2}\kappa h^{-\Delta\alpha-1}\\
\le & 2\eta_{n}^{-4}n\Delta\exp(-c\delta^{2}\eta_{n}^{3}n\Delta h^{2\Delta\alpha})+(n\Delta)^{-1}h^{-\Delta\alpha-1}\to0,
\end{align*}
owing to $n\Delta h^{2\Delta\alpha+1}\to\infty$.
\end{proof}
$\,$
\begin{proof}[Proof of Theorem \ref{thm:rateQuantiles}]
 Proposition~\ref{prop:DistFunct} shows that the numerator in the
error representation (\ref{eq:errorRep}) is of the claimed order.
Moreover, it holds for any $\delta>0$
\[
P(|\hat{\nu}_{h}(\xi^{\pm})-\nu(q_{\tau}^{\pm})|>\delta)\le P(\sup_{|t|<\zeta}|\hat{\nu}_{h}(q_{\tau}^{\pm}+t)-\nu(q_{\tau}^{\pm})|>\delta)+P(|\hat{q}_{\tau,h}^{\pm}-q_{\tau}^{\pm}|\ge\zeta),
\]
where the first term converges to zero by Proposition~\ref{prop:UniformDens}
and the second one tends to zero by Lemma~\ref{lem:consistency}.
Therefore, the denominator in (\ref{eq:errorRep}) can be written
as 
\[
\hat{\nu}_{h}(\xi^{\pm})=\nu(q_{\tau}^{\pm})+o_{P}(1).\qedhere
\]
\end{proof}

\subsection{Proofs for Section \ref{sec:FinancialQuantile}}

To prove Lemma \ref{lem:GaussProc}, we
apply some entropy arguments. To fix the notation, we define
for any (pseudo-)metric $d$ on $\R$ the covering number $N(r,A,d)$
as the smallest number of $d$-balls with radius $r>0$ which is necessary
to cover a subset $A\subset\R$. For $v>0$ the entropy integral is
defined as
\[
J(v,A,d):=\int_{0}^{v}\sqrt{\log N(r,A,d)}\d r,
\]
which is finite for any $v$ if $N(r,A,d)$ grows polynomially in
$r^{-1}$.
\begin{proof}[Proof of Lemma~\ref{lem:GaussProc}]
It\^o's isometry yields
\begin{align}
\Var(\Phi_{n}(u))  =\E[|\Phi_{n}(u)|^{2}]&=n^{-1}u^{2}(u^{2}+1)\E\Big[\Big|\int e^{iux-x}\rho(x)\d W(x)\Big|^{2}\Big]\notag\\
&=n^{-1}u^{2}(u^{2}+1)\|e^{-x}\rho(x)\|_{L^{2}}^{2}\label{eq:VarGaussProc}
\end{align}
and similarly for $\Phi_{n}^{(1)}$ and $\Phi_{n}^{(2)}$ as defined
in (\ref{eq:GaussProc}). From It\^o's isometry and dominated convergence
we conclude that $\Phi_{n}^{(k)},k=1,2,$ are the first and second order
$L^{2}(P)$-derivatives of $\Phi_{n}$. The intrinsic covariance metric
of the Gaussian process $(\int e^{iux-x}x^{k}\rho(x)\d W(x))_{u}$
is given by 
\[
d^{(k)}(u,v):=\E\Big[\Big|\int(e^{iux}-e^{ivx})x^{k}e^{-x}\rho(x)\d W(x)\Big|^{2}\Big]^{1/2}.
\]
Using $\int|x|^{m}e^{-2x}\rho^{2}(x)\d x<\infty$, the entropy integrals
$J(\infty,[-U,U],d^{(k)})$ can be bounded exactly as in the proof
of Proposition~1 by \citet{soehl2014} and are of order $\sqrt{\log U}.$
Dudley's theorem \citep[e.g.][Prop. 3.18]{massart2007} yields then
$\E[\|\Phi_{n}^{(k)}\|_{L^{\infty}[-U,U]}]\lesssim n^{-1/2}U^{2}\sqrt{\log U})$.
\end{proof}

\subsubsection{Convergence rates}

Since the proof strategy is the same for the observation schemes
in Sections~\ref{sec:quantileEst} and \ref{sec:FinancialQuantile},
we will concentrate here on the differences. The estimation error can be decomposed
into $\tilde{N}_{t}(t)-N(t)=B_{n}(t)+\tilde{S}_{n}(t)+V_{n}(t)$ where
$B_{n}(t)$ and $V_{n}(t)$ are given in (\ref{eq:Decomp}) and only
the stochastic error 
\[
\tilde{S}_{n}(t)=\int g_{t}(x)\F^{-1}\Big[\Big(\psi''(u)-\tilde{\psi}_{n}''(u)\Big)\F K(hu)\Big](x)\d x
\]
has a different probabilistic structure. We can apply Lemma~\ref{lem:vola}
and Proposition~\ref{prop:bias} to estimate $B_{n}(t)$ and $V_{n}(t)$. For the sake of brevity, we will frequently write
\[
  (\phi_{T}^{-1}(\tilde{\phi}_{T,n}-\phi_{T}))''
  =\phi_T^{-1}\Phi_n^{(2)}+2(\phi_T^{-1})'\Phi_n^{(1)}+(\phi_T^{-1})''\Phi_n.
\]
This equality is justified in $L^2(P)$-sense, but should merely be understood as notational convention. Linearizing the stochastic error term, we define analogously to (\ref{eq:DecompLin})
\begin{align*}
\tilde{L}_{T,n}(t):= & -\frac{1}{T}\int g_{t}(x){\cal F}^{-1}[\F K(h\bull)(\phi_{T}^{-1}(\tilde{\phi}_{T,n}-\phi_{T}))''](x)\d x.
\end{align*}
Since the $\Phi_{n}^{(k)}$ are almost surely bounded and $\F K$ has
compact support, $\tilde{L}_{T,n}$ is almost surely well defined. The remainder $\tilde S_n-\tilde L_{T,n}$ will be bounded on the event 
\begin{equation}
\Omega_{n,h}:=\Big\{\inf_{u\in I_h}|\tilde{\phi}_{T,n}(u)|\ge n^{-1/2}h^{-2}(\log h^{-1})\Big\}.\label{eq:Omega}
\end{equation}
The order of the remainder in the next lemma corresponds exactly to Lemma~\ref{lem:Remainder} taking the bound from Lemma~\ref{lem:GaussProc} into account.
\begin{lemma}\label{lem:RemainderOptions}
  If $\int|x|^{m}e^{-2x}\rho^{2}(x)\d x<\infty$
  for some $m>4$, then for any sequence $h=h_{n}$ satisfying
  $n^{-1/2}h^{-2}(\log h^{-1})\|\phi_{T}^{-1}\|_{L^{\infty}(I_{h})}\to0$
  as $n\to\infty$ it holds $P(\Omega_{n,h})\to 1$ and uniformly for all L\'evy triplets $(\sigma^2,\gamma,\nu)$
  \begin{align*}
  &\sup_{u\in I_h}\left|\tilde{\psi}_{n}''(u)-\psi''(u)-T^{-1}(\phi_{T}^{-1}(\tilde{\phi}_{T,n}-\phi_{T}))''(u)\right|\\
  & =\mathcal O_P\Big((1+\|\psi'\big\|_{L^{\infty}(I_h)}^{2})n^{-1}h^{-4}\log(h^{-1})\|\phi_{T}^{-1}\|_{L^{\infty}(I_{h})}^{2}\Big).
  \end{align*}
\end{lemma}
\begin{proof}
  $P(\Omega_{n,h})\to 1$ follows immediately from Theorem~1 by \cite{kappusReiss2010}, cf. Lemma 5.1 by \cite{dattnerEtAl2013}. To bound $\tilde{\phi}_{T,n}^{-1}-\phi_{T}^{-1}$, we apply the argument by \citet{neumann1997} and Lemma~\ref{lem:GaussProc}. We obtain on $\Omega_{n,h}$
  \begin{align}
  \|\tilde{\phi}_{T,n}^{-1}-\phi_{T}^{-1}\|_{L^\infty(I_h)}\1_{\Omega_{n,h}}
  &\le\sup_{u\in I_h}\Big(\frac{|\Phi_{n}(u)|}{|\phi_{T}(u)|^{2}}+n^{1/2}h^2|\log h|^{-1}\frac{|\Phi_{n}(u)|^{2}}{|\phi_{T}(u)|^{2}}\Big)\label{eq:NeumannTrick}\\
  &=\mathcal O_P\big(n^{-1/2}h^{-2}|\log h|^{1/2}\|\phi_{T}^{-1}\|^2_{L^{\infty}(I_{h})}\big).\notag
  \end{align}
  A straight forward computation shows on $\Omega_{n,h}$
  \begin{align}
    &T\big(\tilde{\psi}_{n}''-\psi''\big)
    =\frac{\tilde\phi_{T,n}''}{\tilde \phi_{T,n}}-\frac{\phi_T''}{\phi_T}-\Big(\frac{\tilde\phi_{T,n}'}{\tilde\phi_{T,n}}\Big)^2+\Big(\frac{\phi_T'}{\phi_T}\Big)^2\notag\\
    =&\frac{\tilde\phi_{T,n}''-\phi_T''}{\phi_T}
    -\Big(\frac{\tilde\phi_{T,n}'}{\tilde\phi_{T,n}}+\frac{\phi_T'}{\phi_T}\Big)\frac{\tilde\phi_{T,n}'-\phi_T'}{\phi_T}
    +\Big(\frac{\tilde\phi_{T,n}''}{\tilde\phi_{T,n}}-\Big(\frac{\tilde\phi_{T,n}'}{\tilde\phi_{T,n}}+\frac{\phi_T'}{\phi_T}\Big)\frac{\tilde\phi_{T,n}'}{\tilde\phi_{T,n}}\Big)\frac{\phi_T-\tilde\phi_{T,n}}{\phi_T}\notag\\
    =&\frac{\tilde\phi_{T,n}''-\phi_T''}{\phi_T}
    -2\Big(\frac{\phi_{T}'}{\phi_{T}}\Big)\frac{\tilde\phi_{T,n}'-\phi_T'}{\phi_T}
    +\Big(\frac{\phi_{T}''}{\phi_{T}}-2\frac{\phi_T'}{\phi_T}\frac{\phi_{T}'}{\phi_{T}}\Big)\frac{\phi_T-\tilde\phi_{T,n}}{\phi_T}+R_n\notag\\
    =&\phi_T^{-1}\Phi_n^{(2)}+2(\phi_T^{-1})'\Phi_n^{(1)}+(\phi_T^{-1})''\Phi_n+R_n\label{eq:remainderUgly}
  \end{align}
  where $R_n$ is the sum of all second order terms. Using \eqref{eq:NeumannTrick}, formulas \eqref{eq:PhiDerivatives} as well as  Lemma~\ref{lem:GaussProc}, we obtain the claimed order of $R_n$.
\end{proof}

Let us study the linearized stochastic error term. To apply the Lepski method later we need a sharp bound on the variance of $\tilde L_{T,n}$. With the auxiliary functions, $u\in\R$,
\begin{align*}
  \chi_t^{(0)}(u)&:=\F g_{t}(-u)\F K(hu)\big(u(u-i)(\phi_{T}^{-1})''(u)+2(2u-i)(\phi_T^{-1})'(u)+2\phi_T^{-1}(u)\big),\\
  \chi_t^{(1)}(u)&:=\F g_{t}(-u)\F K(hu)\big(2u(iu+1)(\phi_{T}^{-1})'(u)+(4iu+2)\phi_{T}^{-1}(u)\big),\\
  \chi_t^{(2)}(u)&:=u(i-u)\F g_{t}(-u)\F K(hu)\phi_{T}^{-1}(u)
\end{align*}
we define
\begin{align}\label{eq:DefSigma}
  \Sigma_{n,h}(t):=&\frac{1}{2\pi n^{1/2}T}\Big(\|x^{2}e^{-x}\rho(x)\|_{\infty}\|\chi_t^{(2)}\|_{L^2}\notag\\
  &\qquad\quad\qquad+\|xe^{-x}\rho(x)\|_{\infty}\|\chi_t^{(1)}\|_{L^2}+\|e^{-x}\rho(x)\|_{\infty}\|\chi_t^{(0)}\|_{L^2}\Big).
\end{align}
\begin{lemma}\label{lem:LinFin}
  If $\int(1+|x|)^{4}e^{-2x}\rho^{2}(x)\d x<\infty$,
  then $\tilde{L}_{T,h}(t)$ is centered normal. Supposing additionally
  $\|(1\vee x^{2})e^{-x}\rho\|_{\infty}\lesssim1$, it holds
  \begin{align}
  \E[|\tilde{L}_{T,n}(t)|^{2}]^{1/2}
  \le  \Sigma_{n,h}(t)
  \lesssim  n^{-1/2}(t^{-1}\vee t^{-2})\big\|(1+ |u|)|\phi_{T}(u)|^{-1}(1+\psi'(u)^{2})\big\|_{L^{2}(I_{h})}.\nonumber 
  \end{align}
\end{lemma}
\begin{proof}
Since the $\Phi_{n}^{(k)}$ are almost surely bounded, we can apply Plancherel's
identity which yields
\begin{align*}
\E[|\tilde{L}_{T,n}(t)|^{2}] & =\frac{1}{(2\pi T)^2}\E\Big[\Big|\int\F g_{t}(-u)\F K(hu)\phi_{T}^{-1}(u)\Phi_{n}^{(2)}(u)\d u\\
 & \qquad\qquad\quad+2\int\F g_{t}(-u)\F K(hu)(\phi_{T}^{-1})'(u)\Phi_{n}^{(1)}(u)\d u\\
 & \qquad\qquad\quad+\int\F g_{t}(-u)\F K(hu)(\phi_{T}^{-1})''(u)\Phi_{n}(u)\d u\Big|^{2}\Big].
\end{align*}
By continuity and boundedness of $\Phi_{n}^{(k)},k=0,1,2,$ the integral
in $u$ can be approximated with Riemann sums. We conclude first that
$\tilde{L}_{T,h}(t)$ is normally distributed, cf. Section 6.2 in
\citet{soehl2014}, and second that we can exchange the deterministic
integral and the stochastic integral due to the construction of the
Wiener integral as $L^{2}(P)$-limit. Together with It\^o's isometry
and Plancherel's identity we obtain
\begin{align*}
&\E[|\tilde{L}_{T,n}(t)|^{2}]^{1/2}\\
= & \frac{1}{n^{1/2}T}\E\Big[\Big|\int\Big(x^{2}\F^{-1}\chi_t^{(2)}(-x)+x{\cal F}^{-1}\chi_t^{(1)}(-x)+\F^{-1}\chi_t^{(0)}(-x)\Big)e^{-x}\rho(x)\d W(x)\Big|^{2}\Big]^{1/2}\\
= & \frac{1}{n^{1/2}T}\Big(\int\Big|x^{2}\F^{-1}\chi_t^{(2)}(-x)+x{\cal F}^{-1}\chi_t^{(1)}(-x)+\F^{-1}\chi_t^{(0)}(-x)\Big|^{2}e^{-2x}\rho^{2}(x)\d x\Big)^{1/2}\\
\le & \frac{1}{2\pi n^{1/2}T}\Big(\|x^{2}e^{-x}\rho(x)\|_{\infty}\|\chi_t^{(2)}\|_{L^2}+\|xe^{-x}\rho(x)\|_{\infty}\|\chi_t^{(1)}\|_{L^2}+\|e^{-x}\rho(x)\|_{\infty}\|\chi_t^{(0)}\|_{L^2}\Big).
\end{align*}
The formulas (\ref{eq:PhiDerivatives}) and $|\F g_{t}(u)|\lesssim(t^{-1}\vee t^{-2})(1+|u|)^{-1}$ yield the claimed asymptotic bound.
\end{proof}
Now, we can conclude convergence rates for the distribution function
estimator $\tilde{N}_{h}$. 

\begin{prop}
  Suppose $\|(1\vee x^{2})e^{-x}\rho\|_{\infty}\lesssim1$ and $\int(1+|x|)^{m}e^{-2x}\rho^{2}(x)\d x<\infty$
  for some $m>4$. Let $U\subset\R$ be an open set and $\alpha,\beta,s,r,R>0.$
  Let the kernel satisfy \eqref{eq:propKernel} with order $p\ge s+1$. 
  \begin{enumerate}
  \item If $(\sigma^{2},\gamma,\nu)\in\mathcal{D}^{s}(\alpha,2,U,R)$, then
  $|\tilde{N}_{h}(t)-N(t)|=\mathcal{O}_{P}\big(n^{-(s+1)/(2s+2\Delta\alpha+5)}\big)$
  for $h=h_{n}=n^{-1/(2s+2\Delta\alpha+5)}$. 
  \item If $(\sigma^{2},\gamma,\nu)\in\mathcal{E}^{s}(\beta,2,U,r,R)$, then
  $|\tilde{N}_{h}(t)-N(t)|=\mathcal{O}_{P}\big((\log n)^{-(s+1)/\beta}\big)$
  for $h=h_{n}=(\frac{1}{2r}\log n)^{-1/\beta}$. 
  \end{enumerate}
\end{prop}
\begin{proof}
  As in the proof of Proposition \ref{prop:DistFunct} we have
  \begin{align}
  |\tilde{N}_{h}(t)-N(t)|\le & |B_{n}(t)|+|\tilde{S}_{n}(t)|+|V_{n}(t)|\lesssim h^{s+1}+|\tilde{S}_{n}(t)|.\label{eq:ErrDistFunctFin}
  \end{align}
  Lemmas~\ref{lem:RemainderOptions} and \ref{lem:LinFin} yield for
  the stochastic error term
  \begin{align*}
  |\tilde{S}_{n}(t)|& \le|\tilde{L}_{n}(t)|+|\tilde{S}_{n}(t)-\tilde{L}_{n}(t)|\\
  & =\mathcal O_P\Big( n^{-1/2}(1+\|\psi'\big\|_{L^{\infty}(I_h)}^{2})\\
  &\qquad\qquad\times\Big(\big\|(1+|u|)\phi_{T}(u)|^{-1}\big\|_{L^{2}(I_h)}+n^{-1/2}h^{-4}(\log h^{-1})\|\phi_{\Delta}^{-1}\|_{L^{\infty}(I_h)}^{2}\Big)\Big).
  \end{align*}
  In situation (i) we use integrability of $x\nu$ to obtain 
  \begin{align}
  |\tilde{N}_{h}-N|(t)=\mathcal{O}_{P}\Big( & h^{s+1}+n^{-1/2}\big(\|(1+|u|)^{\Delta\alpha+1}\|_{L^{2}(I_h)}+n^{-1/2}(\log h^{-1})h^{-2\Delta\alpha-4}\big)\Big)\notag\\
  =\mathcal{O}_{P}\Big( & h^{s+1}+n^{-1/2}\big(h^{-\Delta\alpha-3/2}+n^{-1/2}(\log h^{-1})h^{-2\Delta\alpha-4}\big)\Big).\label{eq:ErrNmild}
  \end{align}
  Therefore, the optimal bandwidth $h_{n}=n^{-1/(2s+2\Delta\alpha+5)}$
  yields the claimed rate.

  For exponentially decaying characteristic functions in case (ii) we
  infer by $\|\psi'\big\|_{L^{\infty}(I_h)}\lesssim h^{-1}$
  \begin{align}
  |\tilde{N}_{h}-N|(t)= & \mathcal{O}_{P}\Big(h^{s+1}+n^{-1/2}\big(h^{-2}\|(|u|+1)e^{r|u|^{\beta}}\|_{L^{2}(I_h)}+n^{-1/2}h^{-6}(\log h^{-1})e^{rh^{-\beta}}\big)\Big)\notag\\
  = & \mathcal{O}_{P}\Big(h^{s+1}+n^{-1/2}h^{-7/2}(1+n^{-1/2}h^{-5/2}(\log h^{-1}))e^{rh^{-\beta}}\Big),\label{eq:ErrNsev}
  \end{align}
  leading to the rate optimal choice $h_{n}=(\frac{1}{2r}\log n)^{-1/\beta}$.
\end{proof}

Proposition $\ref{prop:UniformDensRates}$ on the density estimator
$\tilde{\nu}_{h}$ is an immediate consequence from the following
proposition.
\begin{prop}\label{prop:UniformDensFin}
  Suppose $\|(1\vee x^{2})e^{-x}\rho\|_{\infty}\lesssim1$ and $\int(1+|x|)^{m}e^{-2x}\rho^{2}(x)\d x<\infty$
  for some $m>4$. Let the kernel satisfy (\ref{eq:propKernel})
  with order $p\ge s$, let $U\subset\R$ be a bounded, open set which
  is bounded away from zero and let $\alpha,\beta,r,R>0$.
  Then we have
  \begin{enumerate}
  \item for any $h\downarrow0$ satisfying $n^{-1/2}h^{-\Delta\alpha-5/2-\delta}\to0$
  for some $\delta>0$ we have uniformly in $(\sigma^{2},\gamma,\nu)\in\mathcal{D}^{s}(\alpha,2,U,R)$ 
  \begin{align*}
  \sup_{t\in U}t^{2}|\tilde{\nu}_{h}(t)-\nu(t)| & =\mathcal O_{P,\mathcal{D}^{s}}\Big( h^{s}+n^{-1/2}(\log n)^{1/2}h^{-\Delta\alpha-5/2}\Big),
  \end{align*}

  \item for any $h\downarrow0$ satisfying $n^{-1/2}e^{(r+\delta)h^{-\beta}}\to0$
  for some $\delta>0$ we have uniformly in $(\sigma^{2},\gamma,\nu)\in\mathcal{E}^{s}(\beta,2,U,r,R)$ 
  \[
  \sup_{t\in U}t^2|\tilde{\nu}_{h}(t)-\nu(t)|=\mathcal O_{P,\mathcal{E}^{s}}\Big( h^{s}+n^{-1/2}(\log n)^{1/2}h^{-9/2}e^{rh^{-\beta}}\Big).
  \]
  \end{enumerate}
\end{prop}
\begin{proof}
As in the proof of Proposition~\ref{prop:UniformDens} we deduce
from Lemma~\ref{lem:RemainderOptions} that
\begin{align*}
\sup_{t\in U}t^2|\tilde{\nu}_{h}(t)-\nu(t)|=&\mathcal O(h^{s})+\sup_{t\in U}t^2|\tilde{L}_{T,n,\nu}(t)|\\
&\qquad+\mathcal O_P\Big((1+\|\psi'\big\|_{L^{\infty}(I_h)}^{2})n^{-1}h^{-5}\log(h^{-1})\|\phi_{\Delta}^{-1}\|_{L^{\infty}(I_h)}^{2}\Big)
\end{align*}
with the linearized stochastic error term
\begin{align*}
\tilde{L}_{T,n,\nu}(t):=-\frac{1}{Tt^2}\F^{-1}\Big[\F K(h\bull)\Big(\frac{\tilde{\phi}_{T,n}-\phi_{T}}{\phi_{T}}\Big)''\Big](t).\\
\end{align*}
We estimate 
\begin{align*}
\E\big[\sup_{t\in U}|t^2\tilde{L}_{T,n,\nu}(t)|\Big] & \le T^{-1}\E\Big[\sup_{t\in U}|{\cal F}^{-1}[m_{T,h}\Phi_n^{(2)}](t)|]+\E[\sup_{t\in U}|{\cal F}^{-1}[m_{T,h}\psi'\Phi_n^{(1)}](t)|]\\
 & \qquad+\E[\sup_{t\in U}|{\cal F}^{-1}[m_{T,h}(\psi''-T(\psi')^2)\Phi_n](t)|\Big]\\
 & =:E_{1}+E_{2}+E_{3}.
\end{align*}
Since all three terms can be estimated analogously, we limit ourselves on
$E_{3}$ which has the largest variance. We estimate uniformly in $t\in U$ with use of Plancherel's
identity, Fubini's theorem and It\^o's isometry
\begin{align}
&\Var\big(\F^{-1}[m_{T,h}(\psi''-T(\psi')^2)\Phi_n](t)\big) \notag\\
 =&\frac{1}{2\pi}\E\Big[\Big|\int e^{-itu}m_{T,h}(u)(\psi''(u)-T\psi'(u)^2)\Phi_n(u)\d u\Big|^{2}\Big]\nonumber \\
 =&n^{-1}\E\Big[\Big|\int\F^{-1}[m_{T,h}(u)(\psi''(u)-T\psi'(u)^2)u(u-i)e^{-itu}](-x)\rho(x)\d W(x)\Big|^{2}\Big]\nonumber \\
 =&n^{-1}\int\big|\F^{-1}[m_{T,h}(u)(\psi''(u)-T\psi'(u)^2)u(u-i)e^{-itu}](-x)\rho(x)\big|^{2}\d x\nonumber \\
 \lesssim& n^{-1}\int(u^{2}+u^{4})|m_{T,h}(u)|^{2}(1+|\psi'(u)|^2)^2\d u=:v(n,h).\label{eq:varDenfin}
\end{align}
Analogously, we can estimate the distance in the intrinsic norm for
any $\delta\in(0,1/2)$ 
\begin{align*}
d(s,t)^{2}:= & \E\Big[\Big|{\cal F}^{-1}[m_{T,h}(\psi''-T(\psi')^2)\Phi_n](t)-{\cal F}^{-1}[m_{T,h}(\psi''-T(\psi')^2)\Phi_n](s)\Big|^{2}\Big]\\
= & n^{-1}\int(1+u^{4})|m_{T,h}(u)|^{2}(1+|\psi'(u)|^2)^2|e^{-itu}-e^{-isu}|^2\d u\\
\lesssim & |t-s|^{2\delta}n^{-1}\int(1+u^{4+2\delta})|m_{T,h}(u)|^{2}(1+|\psi'(u)|^2)^2\d u\\
=: & |t-s|^{2\delta}c_\delta(n,h)
\end{align*}
and thus the covering number is of the order $N(r,U,d)\lesssim(c_\delta(n,h)/r)^{1/(2\delta)}.$
Consequently, the entropy integral can be bounded by
\begin{align}
J\Big(\sqrt{v(n,h)},U,d\Big) & \lesssim\int_{0}^{\sqrt{v(n,h)}}\big(\log c_\delta(n,h)+\log r^{-1}\big)^{1/2}\d r\nonumber \\
 & \lesssim v(n,h)^{1/2}\big(\log c_\delta(n,h)+\log v(n,h)^{-1}\big)^{1/2}.\label{eq:entropyIntDens}
\end{align}
Using this entropy bound, Dudley's theorem \citep[e.g.][Prop. 3.18]{massart2007}
yields 
\[
\E\Big[\sup_{t\in U}|\F^{-1}[m_{T,h}(\psi''-T(\psi')^2)\Phi_n](t)|\Big]\lesssim v(n,h)^{1/2}\big(\log c_\delta(n,h)+\log v(n,h)^{-1}\big)^{1/2}.
\]
Now we can plug in the different assumptions on the decay of $\phi_{T}$
(in particular $\log c_\delta(n,h)$ is smaller than $0$ for $\delta$ sufficiently
small). 
\end{proof}

To finally prove Theorem~\ref{thm:quantileFin}, we can argue as
for Theorem \ref{thm:rateQuantiles}. The only ingredient which remains be shown is uniform convergence $\sup_{|t|>\eta}|\tilde{N}_{h}(t)-N(t)|=o_P(1)$
for the rate optimal bandwidth $h=h_{n}$. Since the remainder $|\tilde{S}_{n}(t)-\tilde{L}_{n}(t)|$
can be bounded uniformly in $t$ using Lemma~\ref{lem:RemainderOptions},
it suffices to show:
\begin{lemma}\label{lem:uniformLtilde}
  If $\|(1\vee x^{2})e^{-x}\rho\|_{\infty}\lesssim1$ and $\int(1+|x|)^{m}e^{-2x}\rho^{2}(x)\d x<\infty$
    for some $m>4$, then it holds uniformly for all $(\sigma^2,\gamma,\nu)$
  \[
  \E[\sup_{|t|\ge\eta}|\tilde{L}_{T,n}(t)|]\lesssim\eta^{-2}n^{-1/2}(\log n)^{1/2}\Big(\int|m_{T,h}(u)|^{2}\big(1+|\psi'(u)|^{2}\big)(1+u^{4})\d u\Big)^{1/2}.
  \]
\end{lemma}
\begin{proof}
  Let us estimate the covering number of $\R\setminus(-\eta,\eta)$
  with respect to the intrinsic metric of the process $\tilde{L}_{T,n}(t)$.
  Similarly to Proposition~\ref{prop:UniformDensFin} we infer 
  \begin{align*}
  d(s,t): & =\E[|\tilde{L}_{T,n}(t)-\tilde{L}_{T,n}(s)|^{2}]^{1/2}\\
  & \lesssim n^{-1/2}\Big(\int|\F[g_{t}-g_{s}](-u)m_{T,h}(u)|^{2}\big(1+|\psi'(u)|^{2}\big)(1+u^{4})\d u\Big)^{1/2}\\
  & \lesssim\|g_{t}-g_{s}\|_{L^{1}}n^{-1/2}\Big(\int|m_{T,h}(u)|^{2}\big(1+|\psi'(u)|^{2}\big)(1+u^{4})\d u\Big)^{1/2}\\
  & =:\|g_{t}-g_{s}\|_{L^{1}}c(n,h).
  \end{align*}
  As in Lemma \ref{lem:consistency}, we see that the covering numbers
  are polynomial: 
  \[
  N\big(r,\R\setminus(-\eta,\eta),d\big)\lesssim\eta_{n}^{-4}c(n,h)^{-2}r^{-2}.
  \]
  With the variance bound $\Sigma_{n,h}:=\sup_{|t|>\eta}\Sigma_{n,h}(t)$ from \eqref{eq:DefSigma} Dudley's theorem and Lemma~\ref{lem:LinFin} yield
  \begin{align*}
  \E[\sup_{|t|\ge\eta}|\tilde{L}_{T,n}(t)|]
  &\lesssim J(\Sigma_{n,h},\R\setminus(-\eta,\eta),d)\\
  &\lesssim\Sigma_{n,h}\big(\log(\eta_{n}^{-2}c(n^{-1/2},h))+\log \Sigma_{n,h}^{-1}\big)^{1/2}.\qedhere
  \end{align*}
\end{proof}

\subsubsection{Adaptive method}

With the previous results at hand the proof  of Theorem~\ref{thm:quantileFinAdapt} is quite similar to the one of Theorem~3.2 by \cite{dattnerEtAl2013} and thus we omit the details. In particular, the bandwidth set $\mathcal{B}_{n}$ fulfills analogous properties as their construction, cf. their Lemma~3.1. Due to $P(\Omega_{n,h})\to1$ for the minimal bandwidth and with $\Omega_{n,h}$ from (\ref{eq:Omega}), it suffices to bound all terms on the complement $\Omega_{n,h}^{c}$. Using (\ref{eq:errorRep}), \eqref{eq:ZEst} and \eqref{eq:ErrDistFunctFin}, the estimation error of $\tilde{q}_{\tau,h}$ can be bounded by 
\begin{align}
|\tilde{q}_{\tau,h}^{\pm}-q_{\tau}^{\pm}|\le & \frac{|B_{n,h}(q_{\tau}^{\pm})|+|\tilde{S}_{n,h}(q_{\tau}^{\pm})|+|V_{n,h}(q_{\tau}^{\pm})|}{|\tilde{\nu}_{n,h}(\xi^{\pm})|}\le\frac{Dh^{s+1}+|\tilde{S}_{n,h}(q_{\tau}^{\pm})|}{|\tilde{\nu}_{n,h}(\xi^{\pm})|}\label{eq:errorDecompTilde}
\end{align}
with probability converging to one and with a deterministic constant $D>0$ involving the bias and the error due to $\sigma^2$. We can verify with use of Proposition~\ref{prop:UniformDensFin}, (\ref{eq:ConsistencyBoundLevy}) and Lemmas~\ref{lem:RemainderOptions}
and \ref{lem:uniformLtilde} for $\eta\in(0,1)$, 
\begin{equation}\label{eq:denominator}
P\Big(\max_{h\in\mathcal{B}_{n}}\sup_{\xi^{\pm}\in[q_{\tau}^{\pm}\wedge\tilde{q}_{\tau,h}^{\pm},q_{\tau}^{\pm}\vee\tilde{q}_{\tau,h}^{\pm}]}|\tilde{\nu}_{n,h}(\xi^{\pm})-\nu(q_{\tau}^{\pm})|>\eta\nu(q_{\tau}^{\pm})\Big)\to0.
\end{equation}
To conclude that $\tilde{V}_{n}^{\pm}(h)$ from \eqref{eq:V} is an appropriate upper bound for $|\tilde{S}_{n,h}(t)|/|\tilde{\nu}_{n,h}(\xi^{\pm})|$
in (\ref{eq:errorDecompTilde}), we have to control $\tilde{S}_{n,h}(q_{\tau}^{\pm})$. We again decompose it into linearization $\tilde{L}_{n,h}(q_{\tau}^{\pm})$ and
remainder. Since $\tilde{L}_{n,h}(q_{\tau}^{\pm})$ is centered and normally distributed
with variance bounded by $\Sigma_{n,h}^{2}(q_\tau^\pm)$ from (\ref{eq:DefSigma}),
the Gaussian concentration and $|\mathcal{B}_{n}|\lesssim\log n$
yield for any $\delta>0$ 
\[
P\Big(\exists h\in\mathcal{B}_{n}:|\tilde{L}_{n,h}(q_{\tau}^{\pm})|>(1+\delta)\sqrt{2\log\log n}\Sigma_{n,h}(q_\tau^\pm)\Big)\to0.
\]
For the remainder, one can show, using (\ref{eq:remainderUgly}) pointwise, $\E[|\tilde{S}_{n,h}(q_\tau^\pm)-\tilde{L}_{n,t}(q_\tau^\pm)|\1_{\Omega_{n,h}}]\lesssim\Sigma_{n,h}(q_\tau^\pm)^{2}h^{-1}$
and thus we conclude for any $\delta>0$
\[
P\Big(\exists h\in\mathcal{B}_{n}:|\tilde{S}_{n,h}(q_\tau^\pm)|>(1+\delta)\sqrt{2\log\log n}\Sigma_{n,h}(q_\tau^\pm)\Big)\to0,
\]
provided that $(\log n)\Sigma_{n,h}(q_\tau^\pm)h^{-1}\to0$. Noting
that the order of $\Sigma_{n,h}(q_\tau^\pm)h^{-1}$ is the same as the order
of the stochastic error of $\tilde{\nu}_{n,h}$, this condition is
satisfied for all $h\in\mathcal{B}_{n}$ by construction. In the next
step, we show that $\Sigma_{n,h}(q_\tau^\pm)$ is reasonably estimated by $\tilde{\Sigma}_{n,h}^\pm$
from (\ref{eq:SigmaTilde}).
\begin{lemma}
In the situation of Theorem~\ref{thm:quantileFinAdapt} we have for any sequence
$h\downarrow0$ satisfying $\inf_{|u|\le 1/h}|\phi_{T}(h)|>n^{-1/2}h^{-2}\log h^{-1}$
that
\[
|\tilde{\Sigma}_{n,h}^\pm-\Sigma_{n,h}(q_\tau^\pm)|=\mathcal O_P\big((\log h^{-1})^{-1}\Sigma_{n,h}(q_\tau^\pm)\big).
\]
\end{lemma}
\begin{proof}
Without loss of generality we only consider $q_\tau^+$. The triangle inequality yields
\begin{align*}
|\tilde{\Sigma}_{n,h}^+-\Sigma_{n,h}(q_\tau^+)|
\le &  \frac{1}{2\pi n^{1/2}T}\Big(\|x^{2}e^{-x}\rho(x)\|_{\infty}\|\tilde\chi_{\tilde q_{\tau,h}^+}^{(2)}-\chi_{q_\tau^+}^{(2)}\|_{L^2}\\
&\quad+\|xe^{-x}\rho(x)\|_{\infty}\|\tilde\chi_{\tilde q_{\tau,h}^+}^{(1)}-\chi_{q_\tau^+}^{(1)}\|_{L^2}+\|e^{-x}\rho(x)\|_{\infty}\|\tilde\chi_{\tilde q_{\tau,h}^+}^{(0)}-\chi_{q_\tau^+}^{(0)}\|_{L^2}\Big).\\
\end{align*}
Hence, it suffices to show
\begin{equation}\label{eq:EstSig}
  \|\tilde\chi_{\tilde q_{\tau,h}^+}^{(k)}-\chi_{q_\tau^+}^{(k)}\|_{L^2}=\mathcal O_P\big((\log h^{-1})^{-1}\|\chi_{q_\tau^+}^{(k)}\|_{L^2}\big),\quad\text{for }k=0,1,2. 
\end{equation}
Let us start with $k=2$ where we have
\begin{align*}
  &\int\big|\tilde\chi_{\tilde q_{\tau,h}^+}^{(2)}(u)-\chi_{q_\tau^+}^{(2)}(u)\big|^2\d u\\
  =&\int(u^2+u^4)|\F K(hu)|^2\big|\F g_{\tilde q_{\tau,h}^+}(-u)\tilde\phi_{T,n}^{-1}(u)-\F g_{q_\tau^+}(-u)\phi_{T}^{-1}(u)\big|^2\d u\\
  \le&\int(u^2+u^4)|\F K(hu)|^2|\phi_T(u)|^{-2}\big|\F g_{\tilde q_{\tau,h}^+}(-u)-\F g_{q_{\tau}^+}(-u)\big|^2\d u\\
  &\quad+\int(u^2+u^4)|\F K(hu)|^2|\F g_{\tilde q_{\tau,h}^+}(-u)|^2\big|\tilde\phi_{T,n}^{-1}(u)-\phi_{T}^{-1}(u)\big|^2\d u\\
  =:&T_1+T_2.
\end{align*}
Using $(g_t-g_s)(x)=x^2\1_{(s,t]}$ for $0<s<t$ and $|\F g_t(u)|\sim(1+|u|)^{-1}$, the first integral can be estimated by
\begin{align*}
 T_1&\le\|g_{\tilde q_{\tau,h}^+}-g_{q_{\tau}^+}\|_{L^1}^2\int(u^2+u^4)|\F K(hu)|^2|\phi_T(u)|^{-2}\d u\\
 &\le\eta_n^{-4}|\tilde q_{\tau,h}^+-q_{\tau}^+|^2h^{-2}\|\chi_{q_\tau^+}^{(2)}\|_{L^2}^2.
\end{align*}
Applying \eqref{eq:ErrNmild}, \eqref{eq:ErrNsev} and \eqref{eq:denominator}, we conclude $T_1=\mathcal O_P((\log h^{-1})^{-2}\|\chi_{q_\tau^+}^{(2)}\|_{L^2}^2)$.

To bound $T_2$, we estimate with \eqref{eq:NeumannTrick} on $\Omega_{n,h}$
\begin{align*}
  T_2\lesssim&\eta_n^{-4}\int(u^2+1)|\F K(hu)|^2|\big|\tilde\phi_{T,n}^{-1}(u)-\phi_{T}^{-1}(u)\big|^2\d u\\
  \le&\eta_n^{-4}\int(u^2+1)|\F K(hu)|^2|\phi_{T}(u)|^{-4}\big(|\Phi_{n}(u)|^{2}+n(u^{2}+u^{4})^{-1}|\Phi_{n}(u)|^{4}\big)\d u
\end{align*}
Using the estimate (\ref{eq:VarGaussProc}), $\Phi_{n}$ is a centered
normal random variable with variance smaller than $n^{-1}u^{2}(u^{2}+1)\|e^{-x}\rho\|_{L^{2}}^{2}$.
Therefore, 
\begin{align*}
T_{2}& =\mathcal O_P\Big(n^{-1}\eta_n^{-4}\int u^2(u^{2}+1)^2|\F K(hu)|^2|\phi_{T}(u)|^{-4}\big)\d u\Big)\\
 & =\mathcal O_P\Big( n^{-1}\eta_n^{-4}h^{-4}\sup_{|v|\le1/h}|\phi_{T}(v)|^{-2}\|\chi_{q_\tau^+}^{(2)}\|_{L^2}^2\Big)\\
 & =\mathcal O_P\Big((\log h^{-1})^{-2}\|\chi_{q_\tau^+}^{(2)}\|_{L^2}^2\Big).
\end{align*}
We conclude \eqref{eq:EstSig} for $k=2$. For $k=0,1$ similar calculations
apply using
\begin{align*}
(\tilde{\phi}_{T,n}^{-1}-\phi_{T}^{-1})' & =\frac{-\Phi_{n}^{(1)}+T\psi'\Phi_{n}}{\tilde{\phi}_{T,n}^{2}}+T\psi'(\phi_{T}^{-1}-\tilde{\phi}_{T,n}^{-1}),\\
(\tilde{\phi}_{T,n}^{-1}-\phi_{T}^{-1})'' & =\frac{-\Phi_{n}^{(2)}+T\psi''\Phi_{n}+T\psi'\Phi_{n}^{(1)}}{\tilde{\phi}_{T,n}^{2}}\\
  &\qquad-\frac{2}{\tilde{\phi}_{T,n}^{2}}\big(1-\frac{\Phi_{n}}{\tilde{\phi}_{T,n}}\big)\big(T\psi'+\frac{\Phi_{n}^{(1)}}{\phi_{T}}\big)\big(-\Phi_{n}^{(1)}+T\psi'\Phi_{n}\big)\\
 & \qquad+T\big(\psi''(\phi_{T}^{-1}-\tilde{\phi}_{T,n}^{-1})+\psi'(\phi_{T}^{-1}-\tilde{\phi}_{T,n}^{-1})'\big).\qedhere
\end{align*}

\end{proof}

With these preparations at hand, we can proceed exactly as in \cite[Sect. 5.2]{dattnerEtAl2013} to see that $\tilde{h}^{\pm}$
mimics the oracle bandwidth which balances the deterministic error
$Dh^{s+1}$ and the stochastic error $\tilde{S}_{n,h}$ in the numerator
of (\ref{eq:errorDecompTilde}). Details are omitted.

\bibliographystyle{apalike}  
\bibliography{library}

\end{document}